\documentclass[12pt, a4paper]{amsart}
\usepackage{amscd,amsmath,amssymb,amsthm,amsfonts}
\usepackage[alphabetic]{amsrefs}
\usepackage{mathrsfs}
\usepackage[shortlabels]{enumitem}
\usepackage[all]{xy}

\usepackage{xcolor}
\usepackage[hypertexnames=true,colorlinks = true, allcolors=blue,linktoc = all, pdffitwindow = false, urlbordercolor = white]{hyperref}%
\usepackage{dbnsymb} 

\usepackage{tikz}
 
\usetikzlibrary{knots}
\usetikzlibrary{decorations.markings}

\usepackage{graphicx} 

 \newcommand{\IN}[0]{\mathbb{N}}

 \newcommand{\IZ}[0]{\mathbb{Z}}

 \renewcommand{\CD}[0]{\mathcal{D}}
 \newcommand{\CF}[0]{\mathcal{F}}

 \newcommand{\CL}[0]{\mathcal{L}}

\newcommand{\tocless}[2]{\bgroup\let\addcontentsline=\nocontentsline#1{#2}\egroup}

\newcommand {\stab}{{\rm Stab}}

\newtheorem{theorem}{Theorem}[section]
\newtheorem*{theorem*}{Theorem}
\newtheorem*{proposition*}{Proposition}
\newtheorem{proposition}[theorem]{Proposition}

\newtheorem*{lemma*}{Lemma}
\newtheorem{example}[theorem]{Example}

\newtheorem{question}{Question}

\newtheorem{remark}[theorem]{Remark}

\newtheorem{convention}{Convention}[section]

\numberwithin{equation}{section}

\begin{document}
\title[An introduction to Thompson knot theory]{An introduction to Thompson knot theory and to Jones subgroups}
\author{Valeriano Aiello} 
\address{Valeriano Aiello,
European Organization for Nuclear Research (CERN) CH-1211 Geneva 23, Switzerland
}\email{valerianoaiello@gmail.com}

\begin{abstract} 
We review a constructions of knots from elements of the Thompson groups due to Vaughan Jones, which comes in two flavours: oriented and unoriented.  
\end{abstract}

\maketitle
\centerline{\it Dedicated to the memory of Vaughan F. R. Jones}\bigskip

\tableofcontents

\section*{Introduction}

%
 %
%

On the Christmas Eve of 2014, the first article \cite{Jo14} of Vaughan Jones on a project centred on the Thompson groups appeared on arXiv and
since then it has been followed by several articles.
The centre of this project was a new powerful machinery that allows one to construct actions of Thompson groups (and more generally of \emph{groups of fractions of certain categories}) 
starting from suitable categories.
There have been developments in several directions. 
A lot of efforts have been devoted to producing  unitary representations of Thompson groups, mainly by means of planar algebras \cite{jo2} and Pythagorean C$^*$-algebras \cite{BJ}, 
see e.g. \cite{Jo14, Jo16, Jo19, ACJ, ABC, AJ, BJ, BJ2, AiCo1, AiCo2, BP, AP}. 
Often these representations are related to notable graph or knot invariants, such as the chromatic polynomial, the Tutte polynomial, the Jones polynomial, the Kauffman bracket, the Homflypt polynomial, to name but a few. 
In another direction, groups can be produced by means of this machinery, see \cite{Brothier1, Brothier2}. 
In this article, we would like to review the developments in two other directions. In particular we would like to talk about the Thompson groups as knot constructors and about Jones subgroups. Our investigations in the latter are mainly motivated by the interest in infinite index maximal subgroups of Thompson groups
and in the classification of the unitary representations introduced by Jones.

In \cite{Jo14} Jones defined a method to produce unoriented knots and links from elements of the Thompson group $F$, which was later extended to the Brown-Thompson group $F_3$.
Since these links do not possess a natural orientation, Jones introduced the oriented subgroups $\vec{F}\leq F$ and $\vec{F}_3\leq F_3$, \cite{Jo14, Jo18}.
For  unoriented links Jones also proved a result analogous to Alexander theorems for braids, that is for every link $L$ there is an element $g$ of the Thompson group whose \emph{closure}
$\CL(g)$ is equal to $L$. In the oriented case a slightly weaker result was proved, namely that the oriented link could be reproduced up to disjoint union with unknots. This result was later strengthened by the author in \cite{A}, where it was shown that every link $\vec{L}$ can be exactly reproduced by choosing a suitable element of the oriented subgroup $\vec{F}$.
With the Thompson group and its oriented subgroup being as good as the braid groups at producing links, it is possible to start a reboot of the theory of braids and links, but with the Thompson group replacing the braid groups.
Like braids, both $F$ and $\vec{F}$ contain interesting \emph{positive} monoids: the monoid of positive words $F_+$ and the monoid of positive oriented words $\vec{F}_+$. The links produced by these monoids were studied by Sebastian Baader and the author 
  in a couple of papers \cite{AB1, AB2}, where it was shown that the links produced by $\vec{F}_+$ are positive 
  (in the sense that all these oriented links admit link diagram where all the crossings are positive)
  and those of $F_+$ arborescent  in the sense of Conway (these knots are also known as algebraic).
 
 Despite being introduced with knot theoretical motivations, $\vec{F}$ and $\vec{F}_3$ turned out to be interesting also from the point of view of group theory.
 In fact, they gave rise to new examples \cite{GS2, TV} of maximal subgroups of infinite index of $F$ and $F_3$, respectively.
 
 We end this introduction with a few words on the structure of the article. Section \ref{sec1} is devoted to introducing the Thompson group and the Brown-Thompson groups. Section \ref{sec2} presents two equivalent methods for producing unoriented links starting from elements of $F$ and $F_3$. Section \ref{sec3} focuses on
positive Thompson links, which are the links produced with elements of $F_+$.
In Section \ref{sec4} the binary and ternary oriented subgroups $\vec{F}$ and $\vec{F}_3$
 are introduced and later in Section 5 they are used to produce oriented links. Section
\ref{sec6} focuses on another Jones’s subgroup: the 3-colorable subgroup $\CF$.

 \section{Preliminaries and notation}\label{sec1}
In this section we recall the definitions of the Thompson group $F$ and of the Brown-Thompson groups $F_k$. The interested reader is referred to \cite{CFP} and \cite{B} for more information on $F$, to \cite{Brown} for $F_k$. 

There are several equivalent definitions of $F$. 
One of them is the following: $F$ is the group of all piecewise linear homeomorphisms of the unit interval $[0,1]$ that are differentiable everywhere except at finitely many dyadic rationals numbers and such that on the intervals of differentiability the derivatives are powers of $2$. We adopt the standard notation: $(f\cdot g)(t)=g(f(t))$.

The Thompson group  $F$ has the following infinite presentation
$$
F=\langle x_0, x_1, x_2, \ldots \; | \; x_nx_k=x_kx_{n+1} \quad \forall \; k<n\rangle\, .
$$
Note that $x_0$ and $x_1$ are enough to generate $F$.
The monoid generated by $x_0, x_1, x_2, \ldots$ is denoted by $F_+$ and its elements are said to be positive. 

Every element $g$ of $F$ can be written in a unique way as 
$$
x_0^{a_0}\cdots x_n^{a_n}x_n^{-b_n}\cdots x_0^{-b_0}
$$
where $a_0$, \ldots , $a_n$, $b_0$, \ldots , $b_n\in \IN_0$,
exactly one between $a_n$ and $b_n$ is non-zero, 
and if $a_i\neq 0$ and $b_i\neq 0$, then $a_{i+1}\neq 0$ or $b_{i+1}\neq 0$ for all $i$. This is the normal form of $g$.

 The projection of $F$ onto its abelianisation is denoted by $\pi: F\to F/[F,F]=\IZ\oplus \IZ$ and it admits a nice interpretation when $F$ is seen as a group of homeomorphisms: $\pi(f)=(\log_2 f'(0),\log_2 f'(1))$. 
 
 A family of groups generalising the Thompson group $F$ are the so-called Brown-Thompson groups.
For any integer $k\geq 2$, the Brown-Thompson group $F_k$ may be defined by the following presentation  
$$
\langle y_0, y_1, \ldots \; | \; y_ny_l=y_ly_{n+k-1} \quad \forall \; l<n\rangle\, .
$$
The elements $y_0, y_1, \ldots , y_{k-1}$ are enough to generate $F_k$.
Note that for $k=2$ we have $F_2=F$.
The monoid generated by $y_0, y_1, y_2, \ldots$ is denoted by $F_{k,+}$
 and it elements are said to be positive. 
In the present article only the monoids $F_+$ and $F_{3,+}$ will play a role.

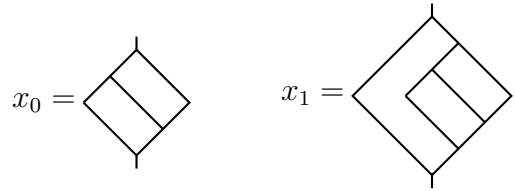
\begin{figure}
\phantom{This text will be invisible} 
\[
\begin{tikzpicture}[x=.35cm, y=.35cm,
    every edge/.style={
        draw,
      postaction={decorate,
                    decoration={markings}
                   }
        }
]

\node at (-1.5,0) {$\scalebox{1}{$x_0=$}$};
\node at (-1.25,-3) {\;};

\draw[thick] (0,0) -- (2,2)--(4,0)--(2,-2)--(0,0);
 \draw[thick] (1,1) -- (2,0)--(3,-1);

 \draw[thick] (2,2)--(2,2.5);

 \draw[thick] (2,-2)--(2,-2.5);

\end{tikzpicture}\qquad
\;\;
\begin{tikzpicture}[x=.35cm, y=.35cm,
    every edge/.style={
        draw,
      postaction={decorate,
                    decoration={markings}
                   }
        }
]

\node at (-3.5,0) {$\scalebox{1}{$x_1=$}$};
\node at (-1.25,-3.25) {\;};

\draw[thick] (2,2)--(1,3)--(-2,0)--(1,-3)--(2,-2);

\draw[thick] (0,0) -- (2,2)--(4,0)--(2,-2)--(0,0);
 \draw[thick] (1,1) -- (2,0)--(3,-1);

 \draw[thick] (1,3)--(1,3.5);
 \draw[thick] (1,-3)--(1,-3.5);

\end{tikzpicture} 
\]
\caption{The generators of $F=F_2$.}\label{genThompsonF}
\end{figure}
\begin{figure}
\phantom{This text will be invisible} 
\[
\begin{tikzpicture}[x=.35cm, y=.35cm,
    every edge/.style={
        draw,
      postaction={decorate,
                    decoration={markings}
                   }
        }
]

\node at (-1.5,0) {$\scalebox{1}{$y_0=$}$};
\node at (-1.25,-3) {\;};

\draw[thick] (0,0) -- (2,2)--(4,0)--(2,-2)--(0,0);
\draw[thick] (1,1) -- (1,0)--(2,-2);
\draw[thick] (1,1) -- (2,0)--(3,-1);
\draw[thick] (2,2) -- (3,0)--(3,-1);

 \draw[thick] (2,2)--(2,2.5);
 \draw[thick] (2,-2)--(2,-2.5);

\end{tikzpicture}
\;\;
\begin{tikzpicture}[x=.35cm, y=.35cm,
    every edge/.style={
        draw,
      postaction={decorate,
                    decoration={markings}
                   }
        }
]

\node at (-1.5,0) {$\scalebox{1}{$y_1=$}$};
\node at (-1.25,-3.25) {\;};

\draw[thick] (0,0) -- (2,2)--(4,0)--(2,-2)--(0,0);
 \draw[thick] (1,0)--(2,-2);
\draw[thick] (3,0)--(2,1) -- (1,0); 
\draw[thick] (2,0)--(3,-1);
\draw[thick] (2,2) -- (2,0);
\draw[thick] (3,0)--(3,-1);

 \draw[thick] (2,2)--(2,2.5);
 \draw[thick] (2,-2)--(2,-2.5);
\end{tikzpicture}
\;\;
\begin{tikzpicture}[x=.35cm, y=.35cm,
    every edge/.style={
        draw,
      postaction={decorate,
                    decoration={markings}
                   }
        }
]

\node at (-3.5,0) {$\scalebox{1}{$y_2=$}$};
\node at (-1.25,-3.25) {\;};

\draw[thick] (2,2)--(1,3)--(-2,0)--(1,-3)--(2,-2);

\draw[thick] (0,0) -- (2,2)--(4,0)--(2,-2)--(0,0);
 \draw[thick] (1,1) -- (2,0)--(3,-1);

\draw[thick] (0,0) -- (2,2)--(4,0)--(2,-2)--(0,0);
\draw[thick] (1,1) -- (1,0)--(2,-2);
\draw[thick] (1,1) -- (2,0)--(3,-1);
\draw[thick] (2,2) -- (3,0)--(3,-1);

\draw[thick] (1,3) -- (-1,0)--(1,-3);


 \draw[thick] (1,3)--(1,3.5);
 \draw[thick] (1,-3)--(1,-3.5);

\end{tikzpicture}
\]
\caption{The generators of $F_3$.}\label{genF3}
\end{figure}
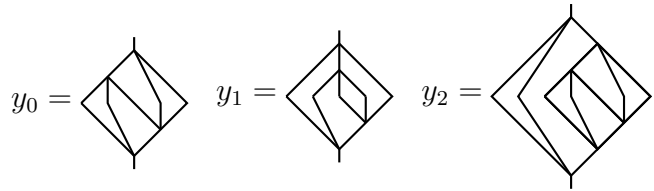

 \begin{figure}
\phantom{This text will be invisible} 
\[
\begin{tikzpicture}[x=1.25cm, y=1.25cm,
    every edge/.style={
        draw,
      postaction={decorate,
                    decoration={markings}
                   }
        }
]

\draw[thick] (0,0)--(.5,.5)--(1,0);
\draw[thick] (0.5,0.5)--(.35,0);
\draw[thick] (0.5,0.5)--(.65,0);
\node at (0,-1.2) {$\;$};
\node at (-.75,0) {$\scalebox{1}{$y_0=$}$};

\draw[thick] (1.25,1.5)--(1.25,1.25);
\draw[thick] (1.25,-1.5)--(1.25,-1.25);

\draw[thick] (0.5,0.5)--(1.25,1.25);
\draw[thick] (1.25,0)--(1.25,1.25);
\draw[thick] (1.75,0)--(1.25,1.25);
\draw[thick] (2.5,0)--(1.25,1.25);

\draw[thick] (2.5,0)--(1.25,-1.25)--(0,0);
\draw[thick] (1.25,-1.25)--(0.35,0);
\draw[thick] (1.25,-1.25)--(0.65,0);

\draw[thick] (1.75,-.75)--(1,0);
\draw[thick] (1.75,-.75)--(1.25,0);
\draw[thick] (1.75,-.75)--(1.75,0);

\end{tikzpicture}
\qquad
\begin{tikzpicture}[x=1.25cm, y=1.25cm,
    every edge/.style={
        draw,
      postaction={decorate,
                    decoration={markings}
                   }
        }
]

\draw[thick] (0.5,0)--(1,.5)--(1.5,0);
\draw[thick] (1,0.5)--(.85,0);
\draw[thick] (1,0.5)--(1.05,0);
\node at (0,-1.2) {$\;$};
\node at (-.75,0) {$\scalebox{1}{$y_1=$}$};

\draw[thick] (0,0)--(1.25,1.25);
\draw[thick] (1,.5)--(1.25,1.25);
\draw[thick] (1.75,0)--(1.25,1.25);
\draw[thick] (2.5,0)--(1.25,1.25);

\draw[thick] (2.5,0)--(1.25,-1.25)--(0,0);
\draw[thick] (1.25,-1.25)--(0.5,0);
\draw[thick] (1.25,-1.25)--(0.85,0);

\draw[thick] (1.75,-.75)--(1.05,0);
\draw[thick] (1.75,-.75)--(1.5,0);
\draw[thick] (1.75,-.75)--(1.75,0);

\draw[thick] (1.25,1.5)--(1.25,1.25);
\draw[thick] (1.25,-1.5)--(1.25,-1.25);

\end{tikzpicture}
\]
\[
\begin{tikzpicture}[x=1.25cm, y=1.25cm,
    every edge/.style={
        draw,
      postaction={decorate,
                    decoration={markings}
                   }
        }
]

\draw[thick] (1,0)--(1.5,.5)--(2,0);
\draw[thick] (1.5,0.5)--(1.35,0);
\draw[thick] (1.5,0.5)--(1.55,0);
\node at (0,-1.2) {$\;$};
\node at (-.75,0) {$\scalebox{1}{$y_2=$}$};

\draw[thick] (0,0)--(1.25,1.25);
\draw[thick] (.5,0)--(1.25,1.25);
\draw[thick] (1.5,.5)--(1.25,1.25);
\draw[thick] (2.5,0)--(1.25,1.25);

\draw[thick] (2.5,0)--(1.25,-1.25)--(0,0);
\draw[thick] (1.25,-1.25)--(0.5,0);
\draw[thick] (1.25,-1.25)--(1,0);

\draw[thick] (1.75,-.75)--(1.35,0);
\draw[thick] (1.75,-.75)--(1.55,0);
\draw[thick] (1.75,-.75)--(2,0);

\node at (1.75,.75) {$\;$};

\draw[thick] (1.25,1.5)--(1.25,1.25);
\draw[thick] (1.25,-1.5)--(1.25,-1.25);

\end{tikzpicture}
\qquad
\begin{tikzpicture}[x=1cm, y=1cm,
    every edge/.style={
        draw,
      postaction={decorate,
                    decoration={markings}
                   }
        }
]

\draw[thick] (.75,1.75)--(.75,2);
\draw[thick] (.75,-1.75)--(.75,-2);

\draw[thick] (0,0)--(.5,.5)--(1,0);
\draw[thick] (0.5,0.5)--(.35,0);
\draw[thick] (0.5,0.5)--(.65,0);
\node at (1.75,-.75) {$\;$};
\node at (-1.75,0) {$\scalebox{1}{$y_3=$}$};

\draw[thick] (0.5,0.5)--(1.25,1.25);
\draw[thick] (1.25,0)--(1.25,1.25);
\draw[thick] (1.75,0)--(1.25,1.25);
\draw[thick] (2.5,0)--(1.25,1.25);

\draw[thick] (-1,0)--(.75,1.75)--(1.25,1.25);
\draw[thick] (-1,0)--(.75,-1.75)--(1.25,-1.25);
\draw[thick] (-.65,0)--(0.75,1.75);
\draw[thick] (-.25,0)--(0.75,1.75);
\draw[thick] (-.65,0)--(0.75,-1.75);
\draw[thick] (-.25,0)--(0.75,-1.75);

\draw[thick] (2.5,0)--(1.25,-1.25)--(0,0);
\draw[thick] (1.25,-1.25)--(0.35,0);
\draw[thick] (1.25,-1.25)--(0.65,0);

\draw[thick] (1.75,-.75)--(1,0);
\draw[thick] (1.75,-.75)--(1.25,0);
\draw[thick] (1.75,-.75)--(1.75,0);

\end{tikzpicture}
\]
\caption{The generators of   $F_4$.}\label{genThompsonF4}
\end{figure}
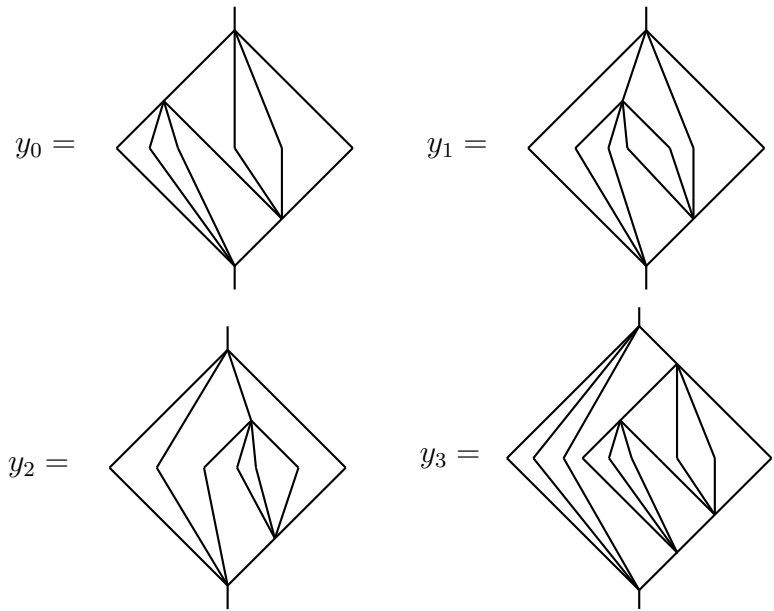

 Going back to $F$,
there is still another description  which is relevant to this paper: the elements of $F$ can be seen as pairs $(T_+,T_-)$ of planar binary rooted trees (with the same number of leaves).  
We draw one tree upside down on top of the other; $T_+$ is  the top tree, while $T_-$ is the bottom tree.
Any pair of binary trees $(T_+,T_-)$ represented in this way is called a binary tree diagram.  Two pairs of binary trees are said to be equivalent if they differ by pairs of opposing carets, namely
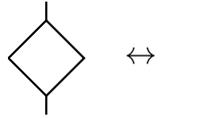
\[\begin{tikzpicture}[x=.5cm, y=.5cm,
    every edge/.style={
        draw,
      postaction={decorate,
                    decoration={markings}
                   }
        }
]

 \draw[thick] (0,0)--(1,1)--(2,0)--(1,-1)--(0,0); 
 \draw[thick] (1,1.5)--(1,1); 
 \draw[thick] (1,-1.5)--(1,-1); 
\node at (0,-1.2) {$\;$};
 \node at (3.5,0) {$\scalebox{1}{$\leftrightarrow$}$};

\end{tikzpicture}
\begin{tikzpicture}[x=.5cm, y=.5cm,
    every edge/.style={
        draw,
      postaction={decorate,
                    decoration={markings}
                   }
        }
]

  \draw[thick] (1,1.5)--(1,-1.5); 
\node at (0,-1.2) {$\;$};
  
\end{tikzpicture}
\]
Every equivalence class of pair of binary trees (i.e. an element of $F$) gives rise to exactly one tree diagram which is reduced, in the sense that 
the number of its vertices is minimal, \cite{B}.  
See Figure \ref{genThompsonF} for the description of $x_0$ and $x_1$ in terms of binary trees.
Thanks to this equivalence relation, the following rule defines the multiplication in    $F$: $(T_+,T)\cdot (T,T_-):=(T_+,T_-)$. The trivial element is represented by any pair $(T,T)$  
and the inverse of $(T_+,T_-)$ is $(T_-,T_+)$. We illustrate how multiplication is performed with $x_0x_1$
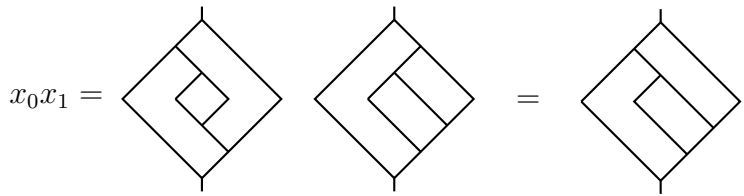
\[
\begin{tikzpicture}[x=.35cm, y=.35cm,
    every edge/.style={
        draw,
      postaction={decorate,
                    decoration={markings}
                   }
        }
]
\node at (-4.5,0) {$\scalebox{1}{$x_0x_1=$}$};
\node at (-1.25,-3.25) {\;};



\draw[thick] (2,2)--(1,3)--(-2,0)--(1,-3)--(2,-2);

\draw[thick] (0,0) -- (1,1);   
\draw[thick] (0,2) -- (1,1);   
\draw[thick] (2,2)--(4,0)--(2,-2)--(0,0);
 \draw[thick] (1,1) -- (2,0)--(1,-1);

 \draw[thick] (1,3)--(1,3.5);
 \draw[thick] (1,-3)--(1,-3.5);

\end{tikzpicture} \quad
\begin{tikzpicture}[x=.35cm, y=.35cm,
    every edge/.style={
        draw,
      postaction={decorate,
                    decoration={markings}
                   }
        }
]

\node at (-1.25,-3.25) {\;};

\draw[thick] (2,2)--(1,3)--(-2,0)--(1,-3)--(2,-2);

\draw[thick] (0,0) -- (2,2)--(4,0)--(2,-2)--(0,0);
 \draw[thick] (1,1) -- (2,0)--(3,-1);

 \draw[thick] (1,3)--(1,3.5);
 \draw[thick] (1,-3)--(1,-3.5);

\end{tikzpicture} 
\quad
\begin{tikzpicture}[x=.35cm, y=.35cm,
    every edge/.style={
        draw,
      postaction={decorate,
                    decoration={markings}
                   }
        }
]

\node at (-2,0) {$\scalebox{1}{$=$}$};
\node at (-1.25,-3) {\;};

\draw[thick] (0,0) -- (3,3)--(6,0)--(3,-3)--(0,0);
 \draw[thick] (3,1) -- (2,0)--(4,-2);

 \draw[thick] (2,2)--(5,-1);

 \draw[thick] (3,-3)--(3,-3.5);
 \draw[thick] (3,3)--(3,3.5);

\end{tikzpicture}
\]

The positive elements of $F$ may always be represented by  a pair of binary trees with  
 bottom tree
having the following shape  
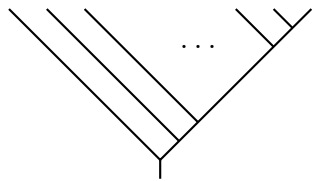
\[
\begin{tikzpicture}[x=.5cm, y=.5cm,
    every edge/.style={
        draw,
      postaction={decorate,
                    decoration={markings}
                   }
        }
]

\draw[thick] (0,0)--(4,-4)--(8,0);
\draw[thick] (4.5,-3.5)--(1,0);
\draw[thick] (5,-3)--(2,0);
\draw[thick] (7,-1)--(6,0);
\draw[thick] (7.5,-.5)--(7,0);
\draw[thick] (4,-4)--(4,-4.5);

\node (aaaa) at (5,-1) {$\scalebox{1}{$\ldots$}$}; 

\end{tikzpicture}
\]

Similarly, the elements of $F_k$ are described by pairs $k$-ary trees (i.e. the trees whose vertices have degree $k+1$, except the leaves having degree $1$).
In this article we will be mainly interested in the Brown-Thompson $F_3$ and $F_4$.
The generators of these groups are displayed in Figures  \ref{genF3} and \ref{genThompsonF4}.
For $F_3$ and $F_4$ pairs of trees are equivalent if they differ by cancellations/additions of pairs of opposing carets 
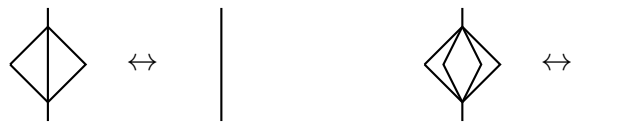
\[
\begin{tikzpicture}[x=.5cm, y=.5cm,
    every edge/.style={
        draw,
      postaction={decorate,
                    decoration={markings}
                   }
        }
]

 \draw[thick] (0,0)--(1,1)--(2,0)--(1,-1)--(0,0); 
 \draw[thick] (1,1.5)--(1,-1.5); 
\node at (0,-1.2) {$\;$};
\node at (3.5,0) {$\scalebox{1}{$\leftrightarrow$}$};

\end{tikzpicture}
\begin{tikzpicture}[x=.5cm, y=.5cm,
    every edge/.style={
        draw,
      postaction={decorate,
                    decoration={markings}
                   }
        }
]

  \draw[thick] (1,1.5)--(1,-1.5); 
\node at (0,-1.2) {$\;$};
  
\end{tikzpicture}
\qquad
\qquad
\qquad
\begin{tikzpicture}[x=.5cm, y=.5cm,
    every edge/.style={
        draw,
      postaction={decorate,
                    decoration={markings}
                   }
        }
]

 \draw[thick] (0,0)--(1,1)--(2,0)--(1,-1)--(0,0); 
 \draw[thick] (1,1.5)--(1,1); 
 \draw[thick] (1,-1.5)--(1,-1);

  \draw[thick] (1,1)--(0.5,0)--(1,-1);
  \draw[thick] (1,1)--(1.5,0)--(1,-1);
 
\node at (0,-1.2) {$\;$};
 \node at (3.5,0) {$\scalebox{1}{$\leftrightarrow$}$};

\end{tikzpicture}
\begin{tikzpicture}[x=.5cm, y=.5cm,
    every edge/.style={
        draw,
      postaction={decorate,
                    decoration={markings}
                   }
        }
]

  \draw[thick] (1,1.5)--(1,-1.5); 
\node at (0,-1.2) {$\;$};
 
\end{tikzpicture}
\]
The bottom ternary tree of a positive element of $F_3$ can be chosen with the following  form 
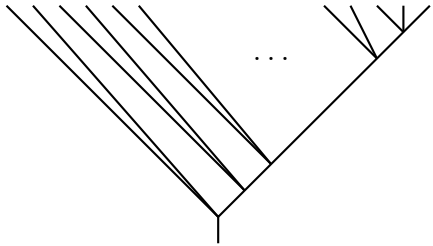
\[
\begin{tikzpicture}[x=.7cm, y=.7cm,
    every edge/.style={
        draw,
      postaction={decorate,
                    decoration={markings}
                   }
        }
]

\draw[thick] (0,0)--(4,-4)--(8,0);
\draw[thick] (4.5,-3.5)--(1,0);
\draw[thick] (5,-3)--(2,0);
\draw[thick] (7,-1)--(6,0);
\draw[thick] (7.5,-.5)--(7,0);
\draw[thick] (4,-4)--(4,-4.5);

\draw[thick] (0.5,0)--(4,-4);
\draw[thick] (1.5,0)--(4.5,-3.5);
\draw[thick] (2.5,0)--(5,-3);
\draw[thick] (7,-1)--(6.5,0);
\draw[thick] (7.5,-.5)--(7.5,0);

\node (aaaa) at (5,-1) {$\scalebox{1}{$\ldots$}$}; 

\end{tikzpicture}
\]

\begin{convention}\label{conventiondrawings}
We draw  $k$-ary trees on the plane with
the roots of our planar $k$-ary trees being drawn as vertices of degree $k+1$. Each $k$-ary tree diagram has the uppermost and lowermost vertices of degree  $1$, which lie respectively on the lines   $y=1$ and $y=-1$.
 The leaves of the trees sit on the $x$-axis, precisely on the non-negative integers. 
\end{convention}
There are an automorphism  and
 two endomorphisms of $F$ that will come in handy later on: 
the flip automorphism and 
the left/right shift homomorphisms  $\sigma, \varphi_L, \varphi_R,: F\to F$.
The flip automorphism $\sigma$ is the order $2$ automorphism  obtained by reflecting tree diagrams about a vertical line, while
the left/right shift homomorphisms $\varphi_L, \varphi_R$ are defined graphically  as
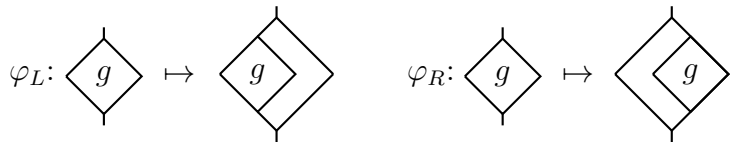
\[
\begin{tikzpicture}[x=1cm, y=1cm,
    every edge/.style={
        draw,
      postaction={decorate,
                    decoration={markings}
                   }
        }
]

\node at (-.45,0) {$\scalebox{1}{$\varphi_L$:}$};

\draw[thick] (0,0)--(.5,.5)--(1,0)--(.5,-.5)--(0,0);
\node at (1.5,0) {$\scalebox{1}{$\mapsto$}$};

\node at (.5,0) {$\scalebox{1}{$g$}$};
\node at (.5,.-.75) {$\scalebox{1}{}$};

 \draw[thick] (.5,.65)--(.5,.5);
 \draw[thick] (.5,-.65)--(.5,-.5);
 
\end{tikzpicture}
\begin{tikzpicture}[x=1cm, y=1cm,
    every edge/.style={
        draw,
      postaction={decorate,
                    decoration={markings}
                   }
        }
]
 
\node at (2.5,0) {$\scalebox{1}{$g$}$};

\draw[thick] (3.5,0)--(2.75,.75)--(2.5,0.5);
\draw[thick] (2.5,-0.5)--(2.75,-.75)--(3.5,0);
\draw[thick] (2,0)--(2.5,.5)--(3,0)--(2.5,-.5)--(2,0);
 
 \node at (1.95,.-.75) {$\scalebox{1}{}$};


 \draw[thick] (2.75,.75)--(2.75,.9);
 \draw[thick] (2.75,-.75)--(2.75,-.9);

\end{tikzpicture}
\qquad
\begin{tikzpicture}[x=1cm, y=1cm,
    every edge/.style={
        draw,
      postaction={decorate,
                    decoration={markings}
                   }
        }
]

\node at (-.45,0) {$\scalebox{1}{$\varphi_R$:}$};

\draw[thick] (0,0)--(.5,.5)--(1,0)--(.5,-.5)--(0,0);
\node at (1.5,0) {$\scalebox{1}{$\mapsto$}$};

\node at (.5,0) {$\scalebox{1}{$g$}$};
\node at (.5,.-.75) {$\scalebox{1}{}$};

 \draw[thick] (.5,.65)--(.5,.5);
 \draw[thick] (.5,-.65)--(.5,-.5);
 
\end{tikzpicture}
\begin{tikzpicture}[x=1cm, y=1cm,
    every edge/.style={
        draw,
      postaction={decorate,
                    decoration={markings}
                   }
        }
]
 
\node at (2.5,0) {$\scalebox{1}{$g$}$};

\draw[thick] (1.5,0)--(2.25,.75)--(3,0)--(2.25,-.75)--(1.5,0);
\draw[thick] (2,0)--(2.5,.5)--(3,0)--(2.5,-.5)--(2,0);
 
 \node at (1.5,.-.75) {$\scalebox{1}{}$};

 \draw[thick] (2.25,.75)--(2.25,.9);
 \draw[thick] (2.25,-.75)--(2.25,-.9);

\end{tikzpicture}
\]
The ranges of $\varphi_L$ and $\varphi_R$ are those elements of $F$ that act trivially on $[1/2,1]$ and $[0,1/2]$, respectively.
Note that $\varphi_R(x_i)=x_{i+1}$ for every $i\in\IN_0$.
Here is $\sigma(x_1)$.
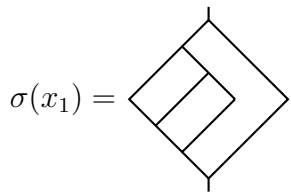
\[
\begin{tikzpicture}[x=.35cm, y=.35cm,
    every edge/.style={
        draw,
      postaction={decorate,
                    decoration={markings}
                   }
        }
]

\node at (-4.5,0) {$\scalebox{1}{$\sigma(x_1)=$}$};
\node at (-1.25,-3.25) {\;};

\draw[thick] (0,2)--(1,3)--(4,0)--(1,-3)--(0,-2);

\draw[thick] (2,0) -- (0,2)--(-2,0)--(0,-2)--(2,0);
\draw[thick] (1,1) -- (-1,-1);

 \draw[thick] (1,3)--(1,3.5);
 \draw[thick] (1,-3)--(1,-3.5);

\end{tikzpicture} 
\]

Some interesting subgroups of $F$ are the so-called rectangular subgroups of $F$. 
They were introduced in \cite{BW} as
\begin{align*}
K_{(a,b)}&:=\{f\in F\; | \; \log_2f'(0)\in a\IZ, \log_2f'(1)\in b\mathbb{Z}\} \qquad a, b\in\IN
\end{align*}
These subgroups can be characterised as the only finite index subgroups isomorphic with $F$ \cite[Theorem 1.1]{BW}.
 
Denote by $W_2$ the set of finite binary words, i.e.  finite sequences of $0$ and $1$.
Let $\IZ[1/2]$ be $\{a/2^k\; |\ ; a, k\in \IZ\}$.
There exists a  map $\rho$ between  
finite binary words 
and
the dyadic rationals in the open unit interval $\CD:=\IZ[1/2]\cap (0,1)$, namely the map $\rho(a_1\ldots a_n):=\sum_{i=1}^n a_i 2^{-i}$ which is bijective when  restricted to finite words ending with $1$ (i.e. $a_n=1$). 
 
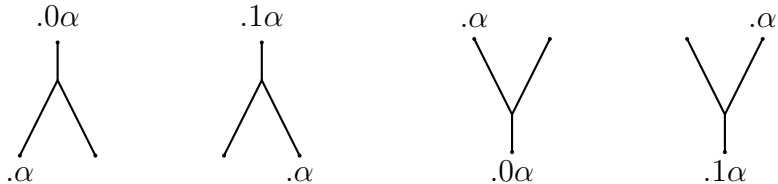
\begin{figure}
\phantom{This text will be invisible} 
 \[
 \begin{tikzpicture}[x=1cm, y=1cm,
    every edge/.style={
        draw,
      postaction={decorate,
                    decoration={markings}
                   }
        }
]

\draw[thick] (0,0) -- (.5,1)--(1,0);
 \draw[thick] (.5,1)--(.5,1.5);

\fill (0,0)  circle[radius=.75pt];
 \fill (1,0)  circle[radius=.75pt];
\fill (.5,1.5)  circle[radius=.75pt];

\node at (.5,1.85) {$\scalebox{1}{$.0\alpha$}$};
\node at (0,-.25) {$\scalebox{1}{$.\alpha$}$};

\end{tikzpicture}
\qquad\qquad
\begin{tikzpicture}[x=1cm, y=1cm,
    every edge/.style={
        draw,
      postaction={decorate,
                    decoration={markings}
                   }
        }
]

\draw[thick] (0,0) -- (.5,1)--(1,0);
 \draw[thick] (.5,1)--(.5,1.5);

\fill (0,0)  circle[radius=.75pt];
 \fill (1,0)  circle[radius=.75pt];
\fill (.5,1.5)  circle[radius=.75pt];

\node at (.5,1.85) {$\scalebox{1}{$.1\alpha$}$};
\node at (1,-.25) {$\scalebox{1}{$.\alpha$}$};

\end{tikzpicture}
\qquad
\qquad
\begin{tikzpicture}[x=1cm, y=1cm,
    every edge/.style={
        draw,
      postaction={decorate,
                    decoration={markings}
                   }
        }
]

\draw[thick] (0,1.5) -- (.5,.5)--(1,1.5);
 \draw[thick] (.5,0)--(.5,.5);

\fill (0,1.5)  circle[radius=.75pt];
\fill (.5,0)  circle[radius=.75pt];
\fill (1,1.5)  circle[radius=.75pt];
 
\node at (.5,-.25) {$\scalebox{1}{$.0\alpha$}$};
\node at (0,1.75) {$\scalebox{1}{$.\alpha$}$};
\node at (1,-.25) {$\scalebox{1}{}$};

\end{tikzpicture}
\qquad\qquad
\begin{tikzpicture}[x=1cm, y=1cm,
    every edge/.style={
        draw,
      postaction={decorate,
                    decoration={markings}
                   }
        }
]

\draw[thick] (0,1.5) -- (.5,.5)--(1,1.5);
 \draw[thick] (.5,0)--(.5,.5);

\fill (0,1.5)  circle[radius=.75pt];
\fill (.5,0)  circle[radius=.75pt];
\fill (1,1.5)  circle[radius=.75pt];
 
\node at (.5,-.25) {$\scalebox{1}{$.1\alpha$}$};
\node at (1,1.75) {$\scalebox{1}{$.\alpha$}$};
\node at (1,-.25) {$\scalebox{1}{}$};

\end{tikzpicture}
\]
\caption{The local rules for computing the action of $F$ on numbers expressed in binary expansion.}\label{compute}
\end{figure}

 The Thompson group $F$ acts by definition on $[0,1]$. 
Now we review this action on the numbers in $[0,1]$
expressed in binary expansion. 
Given a number $t$, it enters into
the top of the top tree in the binary tree diagram, follows a path towards the root of the bottom tree according to the rules portrayed in 
Figure \ref{compute}. What emerges at the bottom is the image of $t$ under the homeomorphism represented by the tree diagram, \cite{BM}. 
Note that there is a change of direction only when the number comes across a vertex of degree $3$ (i.e., the number is unchanged when it comes across a leaf). 

\begin{figure}
\phantom{This text will be invisible} 
 \[
 \begin{tikzpicture}[x=1cm, y=1cm,
    every edge/.style={
        draw,
      postaction={decorate,
                    decoration={markings}
                   }
        }
]

\draw[thick] (0,0) -- (.5,1)--(1,0);
 \draw[thick] (.5,0)--(.5,1.5);

\fill (0,0)  circle[radius=.75pt];
\fill (.5,0)  circle[radius=.75pt];
\fill (1,0)  circle[radius=.75pt];
\fill (.5,1.5)  circle[radius=.75pt];

\node at (.5,1.85) {$\scalebox{1}{$.0\alpha$}$};
\node at (0,-.25) {$\scalebox{1}{$.\alpha$}$};

\end{tikzpicture}
\qquad\qquad
\begin{tikzpicture}[x=1cm, y=1cm,
    every edge/.style={
        draw,
      postaction={decorate,
                    decoration={markings}
                   }
        }
]

\draw[thick] (0,0) -- (.5,1)--(1,0);
 \draw[thick] (.5,0)--(.5,1.5);

\fill (0,0)  circle[radius=.75pt];
\fill (.5,0)  circle[radius=.75pt];
\fill (1,0)  circle[radius=.75pt];
\fill (.5,1.5)  circle[radius=.75pt];

\node at (.5,1.85) {$\scalebox{1}{$.1\alpha$}$};
\node at (.5,-.25) {$\scalebox{1}{$.\alpha$}$};

\end{tikzpicture}
\qquad\qquad
\begin{tikzpicture}[x=1cm, y=1cm,
    every edge/.style={
        draw,
      postaction={decorate,
                    decoration={markings}
                   }
        }
]

\draw[thick] (0,0) -- (.5,1)--(1,0);
 \draw[thick] (.5,0)--(.5,1.5);

\fill (0,0)  circle[radius=.75pt];
\fill (.5,0)  circle[radius=.75pt];
\fill (1,0)  circle[radius=.75pt];
\fill (.5,1.5)  circle[radius=.75pt];

\node at (.5,1.85) {$\scalebox{1}{$.2\alpha$}$};
\node at (1,-.25) {$\scalebox{1}{$.\alpha$}$};

\end{tikzpicture}
 \]
 \[
\begin{tikzpicture}[x=1cm, y=1cm,
    every edge/.style={
        draw,
      postaction={decorate,
                    decoration={markings}
                   }
        }
]

\draw[thick] (0,1.5) -- (.5,.5)--(1,1.5);
 \draw[thick] (.5,0)--(.5,1.5);

\fill (0,1.5)  circle[radius=.75pt];
\fill (.5,0)  circle[radius=.75pt];
\fill (1,1.5)  circle[radius=.75pt];
\fill (.5,1.5)  circle[radius=.75pt];

\node at (.5,-.25) {$\scalebox{1}{$.0\alpha$}$};
\node at (0,1.75) {$\scalebox{1}{$.\alpha$}$};
\node at (1,-.25) {$\scalebox{1}{}$};

\end{tikzpicture}
\qquad\qquad
\begin{tikzpicture}[x=1cm, y=1cm,
    every edge/.style={
        draw,
      postaction={decorate,
                    decoration={markings}
                   }
        }
]

\draw[thick] (0,1.5) -- (.5,.5)--(1,1.5);
 \draw[thick] (.5,0)--(.5,1.5);

\fill (0,1.5)  circle[radius=.75pt];
\fill (.5,0)  circle[radius=.75pt];
\fill (1,1.5)  circle[radius=.75pt];
\fill (.5,1.5)  circle[radius=.75pt];

\node at (.5,-.25) {$\scalebox{1}{$.1\alpha$}$};
\node at (.5,1.75) {$\scalebox{1}{$.\alpha$}$};
\node at (1,-.25) {$\scalebox{1}{}$};

\end{tikzpicture}
\qquad\qquad
\begin{tikzpicture}[x=1cm, y=1cm,
    every edge/.style={
        draw,
      postaction={decorate,
                    decoration={markings}
                   }
        }
]

\draw[thick] (0,1.5) -- (.5,.5)--(1,1.5);
 \draw[thick] (.5,0)--(.5,1.5);

\fill (0,1.5)  circle[radius=.75pt];
\fill (.5,0)  circle[radius=.75pt];
\fill (1,1.5)  circle[radius=.75pt];
\fill (.5,1.5)  circle[radius=.75pt];

\node at (.5,-.25) {$\scalebox{1}{$.2\alpha$}$};
\node at (1,1.75) {$\scalebox{1}{$.\alpha$}$};
\node at (1,-.25) {$\scalebox{1}{}$};

\end{tikzpicture}
\]
\caption{The local rules for computing the action of $F_3$ on numbers expressed in ternary expansion.}\label{compute3}
\end{figure}
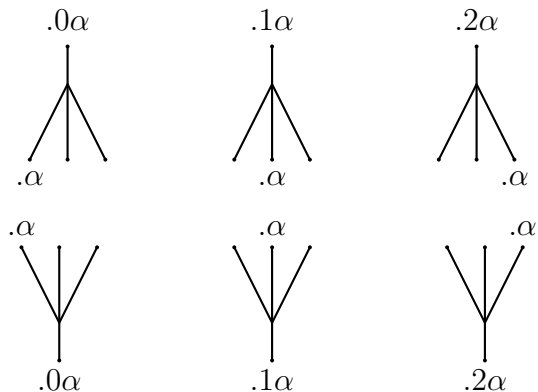

The action of $F_3$ on $[0,1]$ can be describe in a similar way. 
First we express the numbers in ternary expansion (i.e. the digits are only $0$, $1$, $2$). 
Then the number $t$ enters into 
the top of the top tree in the ternary tree diagram, follows a path towards the root of the bottom tree according to the rules portrayed in 
Figure \ref{compute3}.

\section{The construction of knots: the unoriented case}\label{sec2}
Jones introduced two equivalent methods to produce unoriented knots and links from the Thompson groups. Originally these constructions were defined for $F$, \cite{Jo14}, but later they were extended to $F_3$ in \cite{Jo18}.

\tocless\subsection{First method}
We  present this method by taking $x_0x_1$ as  an example.
We will construct a Tait diagram $\Gamma(T_+,T_-)$ from a binary tree diagram $(T_+,T_-)$ in $F$.
Recall that the leaves of $T_+$ and $T_-$ sit on the non-negative integers $\IN_0=\{0,1,2,3,\ldots\}$ on the $x$-axis.
We place the vertices of $\Gamma(T_+,T_-)$ on the half integers, so 
for $x_0x_1$ these points are $(1/2,0)$, $(3/2,0)$, $(5/2,0)$, $(7/2,0)$.
We draw an edge between two of these vertices whenever there is an
edge of the top tree sloping up from left to right (we call them West-North edges, or simply WN$=$\rotatebox[origin=tr]{45}{|}) and 
whenever there is an edge of the bottom tree sloping down from left to right (we refer to them by West-South edges, or just  WS$=$\rotatebox[origin=tr]{-45}{|}).
This is the graph for $x_0x_1$
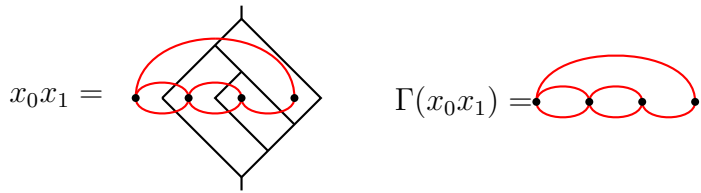
\[
\begin{tikzpicture}[x=.35cm, y=.35cm,
    every edge/.style={
        draw,
      postaction={decorate,
                    decoration={markings}
                   }
        }
]

\node at (-4,0) {$\scalebox{1}{$x_0x_1=$}$};
\node at (-1.25,-3) {\;};

\draw[thick] (0,0) -- (3,3)--(6,0)--(3,-3)--(0,0);
 \draw[thick] (3,1) -- (2,0)--(4,-2);

 \draw[thick] (2,2)--(5,-1);

 \draw[thick] (3,-3)--(3,-3.5);
 \draw[thick] (3,3)--(3,3.5);


\draw[thick, red] (-1,0) to[out=90,in=90] (1,0);
\draw[thick, red] (-1,0) to[out=90,in=180] (2,2.25);
\draw[thick, red] (2,2.25) to[out=0,in=90] (5,0);
\draw[thick, red] (1,0) to[out=90,in=90] (3,0);
\draw[thick, red] (-1,0) to[out=-90,in=-90] (1,0);
\draw[thick, red] (1,0) to[out=-90,in=-90] (3,0);
\draw[thick, red] (3,0) to[out=-90,in=-90] (5,0);

\fill (-1,0)  circle[radius=1.5pt];
\fill (1,0)  circle[radius=1.5pt];
\fill (3,0)  circle[radius=1.5pt];
\fill (5,0)  circle[radius=1.5pt];
\end{tikzpicture}
\qquad
\begin{tikzpicture}[x=.35cm, y=.35cm,
    every edge/.style={
        draw,
      postaction={decorate,
                    decoration={markings}
                   }
        }
] 
 
 \node at (-3.75,0) {$\scalebox{1}{$\Gamma(x_0x_1)=$}$};
 \node at (-2.5,-3) {$\scalebox{1}{$\,$}$};

\draw[thick, red] (-1,0) to[out=90,in=90] (1,0);
\draw[thick, red] (-1,0) to[out=90,in=90] (5,0);
\draw[thick, red] (1,0) to[out=90,in=90] (3,0);
\draw[thick, red] (-1,0) to[out=-90,in=-90] (1,0);
\draw[thick, red] (1,0) to[out=-90,in=-90] (3,0);
\draw[thick, red] (3,0) to[out=-90,in=-90] (5,0);

\fill (-1,0)  circle[radius=1.5pt];
\fill (1,0)  circle[radius=1.5pt];
\fill (3,0)  circle[radius=1.5pt];
\fill (5,0)  circle[radius=1.5pt];

\end{tikzpicture}
\]
There is actually a bijection between the graphs of  the form $\Gamma(T_+,T_-)$ and the pairs of trees $(T_+,T_-)$ of $F$, \cite[Lemma 4.1.4]{Jo14} .
We denote by $\Gamma_+(T_+)$ and $\Gamma_-(T_-)$ the subgraphs of $\Gamma(T_+,T_-)$ contained in the upper and lower-half plane, respectively. 
Since a Tait diagram is a signed graph, we decree that the edges of $\Gamma_+(T_+)$ and $\Gamma_-(T_-)$ are positive and negative, respectively. 


\begin{remark} \label{remarkGamma}
The graphs of the type $\Gamma_\pm(T_\pm)$ may always be assumed to satisfy the following properties
\begin{enumerate}
\item the vertices are $(0,0)$, \ldots , $(N,0)$; 
\item each vertex other than $(0,0)$  is connected to exactly one vertex to its left;
\item each edge $e$ can be parametrized by a function $(x_e(t),y_e(t))$ with $x'_e(t)>0$, for all $t\in [0,1]$, and either $y_e(t)>0$, for all $t\in ]0,1[$ or $y_e(t)<0$, for all $t\in ]0,1[$;
\end{enumerate}
see \cite[Proposition 4.1.3.]{Jo14}. 
 In particular, every vertex (except the leftmost) is the target of exactly two edges, one in the lower half-plane and one in the upper-half plane.
\end{remark}

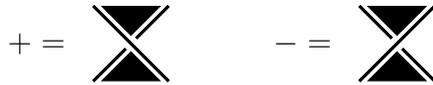
\begin{figure}
 \[\begin{tikzpicture}[every path/.style={very thick}, every node/.style={transform shape, knot crossing, inner sep=1.5pt}]

\node (aaaa) at (-.75,0.5) {$\scalebox{1}{$+=$}$}; 

\node (a1) at (0,0) {}; 
\node (a2) at (0,1) {};
\node (a3) at (1,1) {};
\node (a4) at (1,0) {};
\node (a5) at (0.5,0.5) {};
%
%
\draw (a1.center) .. controls (a1.4 north east) .. (a5) .. controls (a5.4 north east) .. (a3.center);
\draw (a2.center) .. controls (a2.4 south east) .. (a5.center) .. controls (a5.4 south east) .. (a4.center);

  \fill [color=black,opacity=0.2]
               (0.1,0) -- (.5,.4) -- (.9,0) --(.1, 0);

  \fill [color=black,opacity=0.2]
               (0.1,1) -- (.5,.6) -- (.9,1) --(.1, 1); 
\end{tikzpicture}
\qquad 
\qquad\qquad
\begin{tikzpicture}[every path/.style={very thick}, every node/.style={transform shape, knot crossing, inner sep=1.5pt}]

\node (aaaa) at (1.5,0.5) {$\scalebox{1}{$-=$}$}; 

\node (b1) at (2.25,0) {}; 
\node (b2) at (2.25,1) {};
\node (b3) at (3.25,1) {};
\node (b4) at (3.25,0) {};
\node (b5) at (2.75,0.5) {};
 
\draw (b1.center) .. controls (b1.4 north east) .. (b5.center) .. controls (b5.4 north east) .. (b3.center);
\draw (b2.center) .. controls (b2.4 south east) .. (b5) .. controls (b5.4 south east) .. (b4.center);

  \fill [color=black,opacity=0.2]
               (2.35,0) -- (2.75,.4) -- (3.15,0) --(2.35, 0); 

  \fill [color=black,opacity=0.2]
               (2.35,1) -- (2.75,.6) -- (3.15,1) --(2.35, 1);

\end{tikzpicture}\] 
\caption{A positive and a negative crossing. }\label{signscross}
\end{figure}
There are two last steps to be done in order to obtain a link.
First we draw the medial graph $M(\Gamma(T_+,T_-))$ of $\Gamma(T_+,T_-)$. 
In general, given a connected plane graph $G$,
the vertices of its medial graph $M(G)$ sit on every edge of $G$ and an edge of $M(G)$ connects two vertices if they are on adjacent edges of the same face.
Below we will provide an example in our context. 
Now all the vertices of $M(\Gamma(T_+,T_-))$ have degree $4$ and we may make the final step:
turn the vertices into crossings and obtain a link diagram. 
For the vertices in the upper-half plane we use the crossing $\slashoverback$,
while for those in the lower-half plane we use $\backoverslash$. 
We point out that in the 
checkerboard shading of the link diagram obtained with this procedure, 
the crossings corresponding to vertices on edges of $\Gamma_+(T_+)$ are positive and the crossings corresponding to vertices on edges of $\Gamma_-(T_-)$
 are negative (in the sense of Figure \ref{signscross}). 
Here are $M(\Gamma(x_0x_1))$ and $\CL(x_0x_1)$.
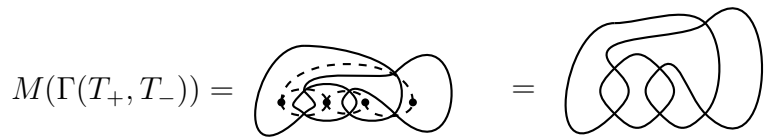
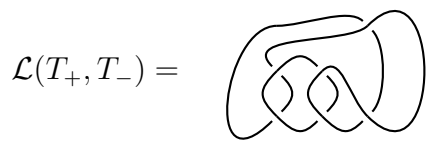
\[
\; \qquad
\qquad
\qquad
\begin{tikzpicture}[x=.6cm, y=.6cm, every path/.style={
 thick}, 
every node/.style={transform shape, knot crossing, inner sep=1.5pt}]
\node (aaaa) at (-3,1) {$\scalebox{1}{$M(\Gamma(T_+,T_-))=$}$}; 

\node (a2) at (1,0.5) {};
\node (a3) at (2,0.5) {};
\node (a4) at (3,0.5) {};
\node (b2) at (1,1) {};
\node (b3) at (2,1) {};
\node (b4) at (3,1.4) {};

\node (appo) at (1,2) {};

\draw (a2.center) .. controls (a2.2 north west) and (b2.4 south west) .. (b2.center);
\draw (a2.center) .. controls (a2.2 north east) and (b2.4 south east) .. (b2.center);

\draw (b2.center) .. controls (b2.4 north) and (b4.2 south west) .. (b4.center);

\draw (b2.center) .. controls (b2.2 north east) and (b3.2 north west) .. (b3.center);
\draw (a2.center) .. controls (a2.2 south east) and (a3.2 south west) .. (a3.center);
\draw (a3.center) .. controls (a3.2 north east) and (b3.2 south east) .. (b3.center);
\draw (a3.center) .. controls (a3.2 north west) and (b3.2 south west) .. (b3.center);
\draw (a2.center) .. controls (a2.16 south west) and (appo.8 west) .. (appo.center);
\draw (appo.center) .. controls (appo.4 east) and (b4.4 north west) .. (b4.center);

\draw (a3.center) .. controls (a3.2 south east) and (a4.2 south west) .. (a4.center);
\draw (a4.center) .. controls (a4.16 south east) and (b4.16 north east) .. (b4.center);
\draw (a4.center) .. controls (a4.2 north east) and (b4.2 south east) .. (b4.center);

\draw (b3.center) .. controls (b3.2 north east) and (a4.2 north west) .. (a4.center);


\draw[thick, dashed] (0.5,.75) to[out=90,in=90] (1.5,.75);
\draw[thick, dashed] (0.5,.75) to[out=-90,in=-90] (1.5,.75);

\draw[thick, dashed] (1.5,.75) to[out=90,in=90] (2.35,.75);
\draw[thick, dashed] (1.5,.75) to[out=-90,in=-90] (2.35,.75);

\draw[thick, dashed] (2.35,.75) to[out=-90,in=-90] (3.4,.75);

\draw[thick, dashed] (0.5,.75) to[out=90,in=90] (3.4,.75);

\fill (0.5,.75) circle[radius=1.5pt];
\fill (1.5,.75)  circle[radius=1.5pt];
\fill (2.35,.75)  circle[radius=1.5pt];
\fill (3.4,.75)  circle[radius=1.5pt];

\end{tikzpicture}
\phantom{This} 
\begin{tikzpicture}[x=.6cm, y=.6cm, every path/.style={
 thick}, 
every node/.style={transform shape, knot crossing, inner sep=1.5pt}]
\node (aaaa) at (-1,0.5) {$\scalebox{1}{$=$}$}; 

\node (a2) at (1,0) {};
\node (a3) at (2,0) {};
\node (a4) at (3,0) {};
\node (b2) at (1,1) {};
\node (b3) at (2,1) {};
\node (b4) at (3,2) {};

\node (appo) at (1,2) {};

\draw (a2.center) .. controls (a2.4 north west) and (b2.4 south west) .. (b2.center);
\draw (a2.center) .. controls (a2.4 north east) and (b2.4 south east) .. (b2.center);

\draw (b2.center) .. controls (b2.8 north west) and (b4.4 south west) .. (b4.center);

\draw (b2.center) .. controls (b2.4 north east) and (b3.4 north west) .. (b3.center);
\draw (a2.center) .. controls (a2.4 south east) and (a3.4 south west) .. (a3.center);
\draw (a3.center) .. controls (a3.4 north east) and (b3.4 south east) .. (b3.center);
\draw (a3.center) .. controls (a3.4 north west) and (b3.4 south west) .. (b3.center);
\draw (a2.center) .. controls (a2.16 south west) and (appo.8 west) .. (appo.center);
\draw (appo.center) .. controls (appo.4 east) and (b4.4 north west) .. (b4.center);

\draw (a3.center) .. controls (a3.4 south east) and (a4.4 south west) .. (a4.center);
\draw (a4.center) .. controls (a4.16 south east) and (b4.16 north east) .. (b4.center);
\draw (a4.center) .. controls (a4.4 north east) and (b4.4 south east) .. (b4.center);

\draw (b3.center) .. controls (b3.4 north east) and (a4.4 north west) .. (a4.center);

\end{tikzpicture}
\phantom{Ths} \]
\[
\begin{tikzpicture}[x=.6cm, y=.6cm, every path/.style={
 thick}, 
every node/.style={transform shape, knot crossing, inner sep=1.5pt}]

\node (aaaa) at (-3,1) {$\scalebox{1}{$\CL(T_+,T_-)=$}$}; 

\node (a2) at (1,0) {};
\node (a3) at (2,0) {};
\node (a4) at (3,0) {};
\node (b2) at (1,1) {};
\node (b3) at (2,1) {};
\node (b4) at (3,2) {};

\node (appo) at (1,2) {};

\draw (a2.center) .. controls (a2.4 north west) and (b2.4 south west) .. (b2.center);
\draw (a2) .. controls (a2.4 north east) and (b2.4 south east) .. (b2);

\draw (b2) .. controls (b2.8 north west) and (b4.4 south west) .. (b4.center);

\draw (b2.center) .. controls (b2.4 north east) and (b3.4 north west) .. (b3);
\draw (a2.center) .. controls (a2.4 south east) and (a3.4 south west) .. (a3);
\draw (a3) .. controls (a3.4 north east) and (b3.4 south east) .. (b3);
\draw (a3.center) .. controls (a3.4 north west) and (b3.4 south west) .. (b3.center);
\draw (a2) .. controls (a2.16 south west) and (appo.8 west) .. (appo.center);
\draw (appo.center) .. controls (appo.4 east) and (b4.4 north west) .. (b4);

\draw (a3.center) .. controls (a3.4 south east) and (a4.4 south west) .. (a4);
\draw (a4.center) .. controls (a4.16 south east) and (b4.16 north east) .. (b4.center);
\draw (a4) .. controls (a4.4 north east) and (b4.4 south east) .. (b4);

\draw (b3.center) .. controls (b3.4 north east) and (a4.4 north west) .. (a4.center);

\end{tikzpicture}
\]

\tocless\subsection{Second method}
In this section we describe an equivalent procedure to obtain links from elements of $F$,  \cite{Jo14, Jo18}.
The advantage of this description is that it can be readily extended to $F_3$.
We start with a binary tree diagram in $F$. 
The first operation is to turn all the $3$-valent vertices into $4$-valent by adding additional edges below each vertex of degree $3$ in the top tree and above each vertex of degree $3$ of the bottom tree, which we join in the only planar possible way.
The second operation is to draw an edge between the two roots of the trees. 
The third and last operation is to turn all the $4$-valent vertices into crossings as shown in Figure \ref{rulesknot}: the vertices and the four incident edges are replaced by "forks", 
see leftmost illustration of Fig. \ref{rulesknot}.
We exemplify this procedure with $x_0x_1$. 
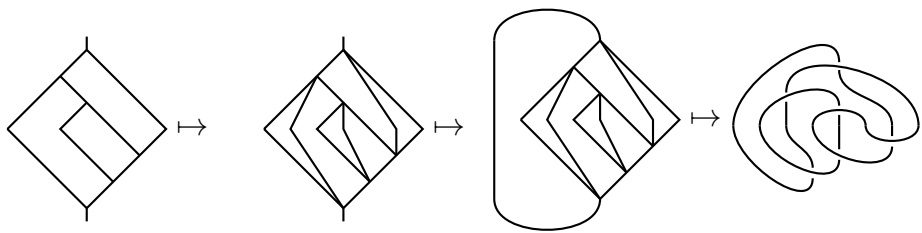
\begin{figure}
\[
\begin{tikzpicture}[x=.3cm, y=.3cm,
    every edge/.style={
        draw,
      postaction={decorate,
                    decoration={markings}
                   }
        }
]

\draw[thick] (1,1)--(1,2);
\draw[thick] (0,0) --(1,1)--(2,0);

\node at (0,-1.2) {$\;$};

\end{tikzpicture}\quad
\begin{tikzpicture}[x=.3cm, y=.3cm,
    every edge/.style={
        draw,
      postaction={decorate,
                    decoration={markings}
                   }
        }
]

\draw[thick] (1,0)--(1,2);
\draw[thick] (0,0) --(1,1)--(2,0);

\node at (0,-1.2) {$\;$};
\node at (-2,1) {$\scalebox{1}{$\mapsto\; $}$};

\end{tikzpicture}\qquad
\begin{tikzpicture}[x=.3cm, y=.3cm,
    every edge/.style={
        draw,
      postaction={decorate,
                    decoration={markings}
                   }
        }
]

\draw[thick] (1,0)--(1,2);
\draw[thick] (0,0) --(1,1)--(2,0);

\node at (0,-1.2) {$\;$};

\end{tikzpicture}\quad
\begin{tikzpicture}[x=.3cm, y=.3cm,
    every edge/.style={
        draw,
      postaction={decorate,
                    decoration={markings}
                   }
        }
]

\draw[thick] (1,0)--(1,.35);
\draw[thick] (1,.75)--(1,2);
\draw[thick] (0,0) to[out=90,in=90] (2,0);

\node at (-2,1) {$\scalebox{1}{$\mapsto\; $}$};

\node at (0,-1.2) {$\;$};

\end{tikzpicture}
\qquad \begin{tikzpicture}[x=.3cm, y=.3cm,
    every edge/.style={
        draw,
      postaction={decorate,
                    decoration={markings}
                   }
        }
]

\draw[thick] (1,0)--(1,2);
\draw[thick] (0,2) --(1,1)--(2,2);

\node at (0,-1.2) {$\;$};

\end{tikzpicture}\quad
\begin{tikzpicture}[x=.3cm, y=.3cm,
    every edge/.style={
        draw,
      postaction={decorate,
                    decoration={markings}
                   }
        }
]

\draw[thick] (1,2)--(1,1.65);
\draw[thick] (1,1.25)--(1,0);
\draw[thick] (0,2) to[out=-90,in=-90] (2,2);

\node at (-2,1) {$\scalebox{1}{$\mapsto\; $}$};

\node at (0,-1.2) {$\;$};

\end{tikzpicture}
\]
\caption{The rules needed for obtaining $\CL(g)$. }\label{rulesknot}
\end{figure}
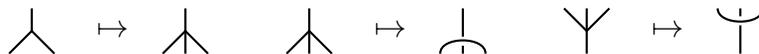
\[
\begin{tikzpicture}[x=.35cm, y=.35cm,
    every edge/.style={
        draw,
      postaction={decorate,
                    decoration={markings}
                   }
        }
]
 
\draw[thick] (0,0) -- (3,3)--(6,0)--(3,-3)--(0,0);
 \draw[thick] (3,1) -- (2,0)--(4,-2);

 \draw[thick] (2,2)--(5,-1);

 \draw[thick] (3,-3)--(3,-3.5);
 \draw[thick] (3,3)--(3,3.5);

\node at (7,0) {$\scalebox{1}{$\mapsto$}$};

\node at (0,-1.2) {\phantom{$\frac{T_+}{T_-}=$}};
\node at (0,-3) {\phantom{$\frac{T_+}{T_-}=$}};

\end{tikzpicture}
\begin{tikzpicture}[x=.35cm, y=.35cm,
    every edge/.style={
        draw,
      postaction={decorate,
                    decoration={markings}
                   }
        }
]

 \draw[thick] (3,3.5) -- (3,3);  
 \draw[thick] (3,-3.5) -- (3,-3);  

\draw[thick] (0,0) -- (3,3)--(6,0)--(3,-3)--(0,0);
 \draw[thick] (3,1) -- (2,0)--(4,-2);

 \draw[thick] (2,2)--(5,-1);

\node at (7,0) {$\scalebox{1}{$\mapsto$}$};

 \draw[thick] (2,2)--(1,0)--(3,-3);
 \draw[thick] (3,1)--(3,0)--(4,-2);
 \draw[thick] (5,-1)--(5,0)--(3,3);

\node at (0,-1.2) {\phantom{$\frac{T_+}{T_-}=$}};
\node at (0,-3) {\phantom{$\frac{T_+}{T_-}=$}};

\end{tikzpicture}
\begin{tikzpicture}[x=.35cm, y=.35cm,
    every edge/.style={
        draw,
      postaction={decorate,
                    decoration={markings}
                   }
        }
]

 \draw[thick] (-1,3) to[out=90,in=90] (3,3);  
 \draw[thick] (-1,-3) to[out=-90,in=-90] (3,-3);  
\draw[thick] (-1,-3) --(-1,3);

\draw[thick] (0,0) -- (3,3)--(6,0)--(3,-3)--(0,0);
 \draw[thick] (3,1) -- (2,0)--(4,-2);

 \draw[thick] (2,2)--(5,-1);

\node at (7,0) {$\scalebox{1}{$\mapsto$}$};

 \draw[thick] (2,2)--(1,0)--(3,-3);
 \draw[thick] (3,1)--(3,0)--(4,-2);
 \draw[thick] (5,-1)--(5,0)--(3,3);

\node at (0,-1.2) {\phantom{$\frac{T_+}{T_-}=$}};
\node at (0,-3) {\phantom{$\frac{T_+}{T_-}=$}};

\end{tikzpicture}
\begin{tikzpicture}[x=.35cm, y=.35cm,
    every edge/.style={
        draw,
      postaction={decorate,
                    decoration={markings}
                   }
        }
]

\draw[thick] (2,0) to[out=90,in=90] (4,0);   
\draw[thick] (0,0) to[out=90,in=90] (3,.75);   
\draw[thick] (1,1.25) to[out=90,in=90] (6,0);   
\draw[thick] (3,2.4) to[out=90,in=90] (-1,0);   
\draw[thick] (1,.9)--(1,0);
\draw[thick] (3,.4)--(3,0);
\draw[thick] (3,2.1) to[out=-90,in=90] (5,0);

\draw[thick] (4,0) to[out=-90,in=-90] (6,0);   
\draw[thick] (2,0) to[out=-90,in=-90] (5,-.8);   
\draw[thick] (0,0) to[out=-90,in=-90] (3,-1.3);   
\draw[thick] (-1,0) to[out=-90,in=-90] (2,-2);   
\draw[thick] (5,0)--(5,-.4);
\draw[thick] (3,0)--(3,-.9);
\draw[thick] (1,0) to[out=-90,in=90] (2,-1.6);
 
\node at (0,-4) {$\scalebox{1}{$\,$}$};

\node at (0,-1.2) {$\;$};

\end{tikzpicture}
\qquad
\]
\begin{remark}
In the first step 
we are actually using an injective group homomorphism $\iota: F\to F_3$, see Figure \ref{iota}. This map takes a binary tree diagram, 
and returns $4$-valent tree diagrams.
\begin{figure}
\phantom{This text will be invisible} 
\[\begin{tikzpicture}[x=.5cm, y=.5cm,
    every edge/.style={
        draw,
      postaction={decorate,
                    decoration={markings}
                   }
        }
]

 \draw[thick] (0,0)--(1,1)--(2,0); 
 \draw[thick] (1,1.5)--(1,1); 

\node at (3.5,.5) {$\scalebox{1}{$\mapsto$}$};

\end{tikzpicture}\; \;
\begin{tikzpicture}[x=.5cm, y=.5cm,
    every edge/.style={
        draw,
      postaction={decorate,
                    decoration={markings}
                   }
        }
]

 \draw[thick] (0,0)--(1,1)--(2,0); 
 \draw[thick] (1,1.5)--(1,0); 
 
\end{tikzpicture}
\]

\[
\begin{tikzpicture}[x=.35cm, y=.35cm,
    every edge/.style={
        draw,
      postaction={decorate,
                    decoration={markings}
                   }
        }
]

\node at (-1.25,-3.25) {\;};

\draw[thick] (2,2)--(1,3)--(-2,0)--(1,-3)--(2,-2);

\draw[thick] (0,0) -- (2,2)--(4,0)--(2,-2)--(0,0);
 \draw[thick] (1,1) -- (2,0)--(3,-1);

 \draw[thick] (1,3)--(1,3.5);
 \draw[thick] (1,-3)--(1,-3.5);

\end{tikzpicture} \; \;
\begin{tikzpicture}[x=.35cm, y=.35cm,
    every edge/.style={
        draw,
      postaction={decorate,
                    decoration={markings}
                   }
        }
]

\node at (-3.5,0) {$\scalebox{1}{$\mapsto$}$};
\node at (-1.25,-3.25) {\;};

\draw[thick] (2,2)--(1,3)--(-2,0)--(1,-3)--(2,-2);

\draw[thick] (0,0) -- (2,2)--(4,0)--(2,-2)--(0,0);
 \draw[thick] (1,1) -- (2,0)--(3,-1);

\draw[thick] (0,0) -- (2,2)--(4,0)--(2,-2)--(0,0);
\draw[thick] (1,1) -- (1,0)--(2,-2);
\draw[thick] (1,1) -- (2,0)--(3,-1);
\draw[thick] (2,2) -- (3,0)--(3,-1);

\draw[thick] (1,3) -- (-1,0)--(1,-3);

 \draw[thick] (1,3)--(1,3.5);
 \draw[thick] (1,-3)--(1,-3.5);

\end{tikzpicture}
\]
\caption{The map $\iota: F\to F_3$.}\label{iota}
\end{figure}
This map was originally defined   
by Jones in \cite[Section 4]{Jo18}. Therefore, the construction of knots can be extended to $F_3$ just by skipping the first step.
\end{remark}
We draw the Tait graph for the links corresponding to $F_3$ in the plane, the vertices sitting on the $x$-axis, half of the edges in the upper-half-plane, the other half in the lower-half-plane.
When restricted to $F$,
the Tait diagram of the link diagram obtained in this way is exactly the graph $\Gamma(T_+,T_-)$ described in the previous section.
One of the differences between the Tait graphs of the elements $F$ and those of $F_3$ is that for the elements of$F_3$ is that the edges in the lower-half-plane (upper-half-plane, respectively) are not necessarily negative (positive, respectively). 

\begin{example}[The $4_1$ knot]
Consider the element $g=x_0x_2^2x_5x_6(x_4x_6x_7)^{-1} \in F$ whose image $\iota(g)=(T_+,T_-)=y_0y_4^2y_{10}y_{12}(y_8y_{12}y_{14})^{-1}\in F_3$ is described by the following pair of ternary trees
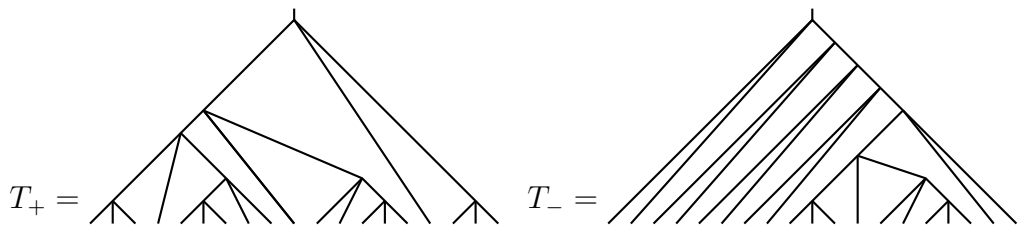
\[\begin{tikzpicture}[x=.3cm, y=.3cm,
    every edge/.style={
        draw,
      postaction={decorate,
                    decoration={markings}
                   }
        }
]

\draw[thick] (10,9)--(10,9.5);
\draw[thick] (1,0)--(10,9)--(19,0);
\draw[thick] (2,0)--(2,1);
\draw[thick] (2,1)--(3,0);
\draw[thick] (4,0)--(5,4);
\draw[thick] (5,4)--(7,2)--(9,0);
\draw[thick] (5,0)--(7,2);
\draw[thick] (6,0)--(6,1);
\draw[thick] (7,0)--(6,1);
\draw[thick] (8,0)--(7,2);
\draw[thick] (10,0)--(6,5);
\draw[thick] (15,0)--(14,1)--(13,2)--(6,5);
\draw[thick] (10,0)--(6,5);
\draw[thick] (11,0)--(13,2);
\draw[thick] (12,0)--(13,2);
\draw[thick] (13,0)--(14,1);
\draw[thick] (14,0)--(14,1);
\draw[thick] (16,0)--(10,9);
\draw[thick] (17,0)--(18,1);
\draw[thick] (18,0)--(18,1);

\node at (-1,1) {$\scalebox{1}{$T_+=$}$};

\node at (0,-1.2) {$\;$};

\end{tikzpicture}\; \; 
\begin{tikzpicture}[x=.3cm, y=.3cm,
    every edge/.style={
        draw,
      postaction={decorate,
                    decoration={markings}
                   }
        }
]

\draw[thick] (10,9)--(10,9.5);
\draw[thick] (1,0)--(10,9)--(19,0);
\draw[thick] (2,0)--(10,9);
\draw[thick] (3,0)--(11,8);
\draw[thick] (4,0)--(11,8);
\draw[thick] (5,0)--(12,7);
\draw[thick] (6,0)--(12,7);
\draw[thick] (7,0)--(13,6);
\draw[thick] (8,0)--(13,6);
\draw[thick] (9,0)--(14,5);
\draw[thick] (10,0)--(10,1);
\draw[thick] (11,0)--(10,1);
\draw[thick] (12,0)--(12,3);
\draw[thick] (17,0)--(15,2)--(12,3);
\draw[thick] (13,0)--(15,2);
\draw[thick] (14,0)--(15,2);
\draw[thick] (15,0)--(16,1);
\draw[thick] (16,0)--(16,1);
\draw[thick] (18,0)--(14,5);

\node at (-1,1) {$\scalebox{1}{$T_-=$}$};

\node at (0,-1.2) {$\;$};

\end{tikzpicture}
\] 
After applying some Reidemeister moves (to be precise a sequence consisting of five Reidemeister moves
of type II and four of type I), one sees that $\CL(T_+,T_-)$ is the $4_1$ knot.
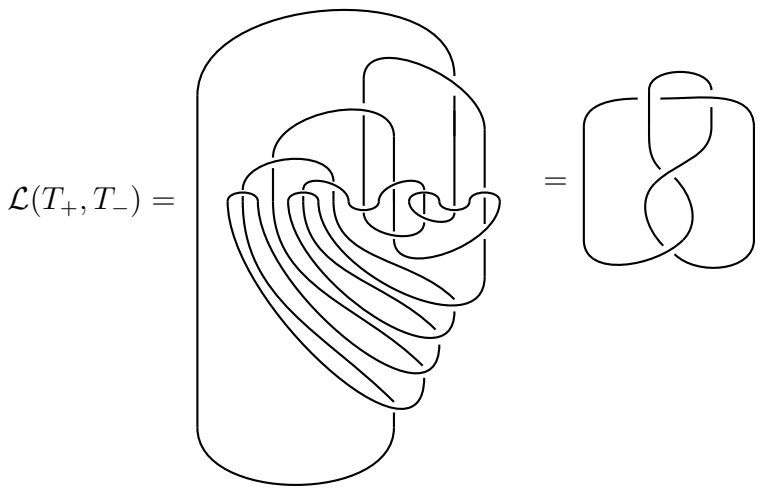
\[\begin{tikzpicture}[x=.2cm, y=.2cm,
    every edge/.style={
        draw,
      postaction={decorate,
                    decoration={markings}
                   }
        }
]

\draw[thick] (-1,-15) to[out=90,in=-90] (-1,7); 
\draw[thick] (-1,7) to[out=90,in=90] (16,8.3); 
\draw[thick] (1,0) to[out=90,in=90] (3,0); 
\draw[thick] (2,0) to[out=90,in=-90] (2,0.35); 
\draw[thick] (2,.75) to[out=90,in=90] (8,1.3); 
\draw[thick] (4,0) to[out=90,in=90] (4,2); 
\draw[thick] (4,2.9) to[out=90,in=90] (12,4.3); 
\draw[thick] (12,1.3) to[out=90,in=-90] (12,4.3); 
\draw[thick] (5,0) to[out=90,in=90] (7,0); 
\draw[thick] (6,0.75) to[out=90,in=90] (9,0); 
\draw[thick] (6,0.35) to[out=-90,in=90] (6,0); 
\draw[thick] (8,0) to[out=90,in=90] (8,0.6); 
\draw[thick] (10,-0.3) to[out=90,in=-90] (10,5.7); 
\draw[thick] (10,6.5) to[out=90,in=-90] (10,8.1); 
\draw[thick] (10,8.1) to[out=90,in=90] (18,4.75); 
\draw[thick] (11,0) to[out=90,in=90] (14,.75); 
\draw[thick] (12,0) to[out=90,in=90] (12,.75); 
\draw[thick] (12,-2) to[out=90,in=-90] (12,.75); 
\draw[thick] (14,-0.85) to[out=90,in=90] (14,.35); 
\draw[thick] (13,0) to[out=90,in=90] (15,0); 
\draw[thick] (16,-0.3) to[out=90,in=-90] (16,6); 
\draw[thick] (16,4.3) to[out=90,in=-90] (16,7); 
\draw[thick] (17,0) to[out=90,in=90] (19,0); 
\draw[thick] (18,-1.5) to[out=90,in=-90] (18,0.25); 
\draw[thick] (18,-2.5) to[out=-90,in=90] (18,-5); 
\draw[thick] (18,.9) to[out=90,in=-90] (18,4.75);

\node at (-8,0) {$\scalebox{1}{$\CL(T_+,T_-)=$}$};

\draw[thick] (1,0) to[out=-90,in=-90] (14,-11.7);
\draw[thick] (2,0) to[out=-90,in=135] (12,-13.3);  
\draw[thick] (3,0) to[out=-90,in=-90] (15,-9.5); 
\draw[thick] (4,0) to[out=-90,in=135] (13.9,-10.9); 
\draw[thick] (5,0) to[out=-90,in=-90] (16,-7.3); 
\draw[thick] (6,0) to[out=-90,in=135] (14.75,-8.5); 
\draw[thick] (7,0) to[out=-90,in=-90] (18,-5); 
\draw[thick] (8,0) to[out=-90,in=135] (16,-6.5); 
\draw[thick] (9,0) to[out=-90,in=-90] (11,0); 
\draw[thick] (10,-0.75) to[out=-90,in=-90] (14,-1.4); 
\draw[thick] (12,-2.5) to[out=-90,in=-90] (19,0); 
\draw[thick] (12,-14) to[out=-90,in=90] (12,-15);  
\draw[thick] (13,0) to[out=-90,in=-90] (16,-0.75); 
\draw[thick] (15,0) to[out=-90,in=-90] (17,0); 
\draw[thick] (-1,-15) to[out=-90,in=-90] (12,-15);

\node at (0,-1.2) {$\;$};

\end{tikzpicture} 
\quad
\begin{tikzpicture}[x=.75cm, y=.75cm,
    every edge/.style={
        draw,
      postaction={decorate,
                    decoration={markings}
                   }
        }
]
 
\draw[thick] (-.25,0) to[out=-135,in=135] (-.1,-1.1);  
\draw[thick] (-.25,0) to[out=45,in=-90] (.75,1);  

\draw[thick] (.1,.1) to[out=-45,in=45] (.25,-1);  
\draw[thick] (.25,-1) to[out=-135,in=-90] (-1.5,-1);  

\draw[thick] (.1,-1.25) to[out=-45,in=-90] (1.5,-1);  
\draw[thick] (1.5,-1)--(1.5,1); 

\draw[thick] (-1.5,-1)--(-1.5,1); 

\draw[thick] (.75,1)--(.75,1.35); 

\draw[thick] (-.15,.25) to[out=135,in=-90] (-.35,1); 

\draw[thick] (.75,1.65) to[out=90,in=90] (-.35,1.65); 

\draw[thick] (-1.5,1) to[out=90,in=180] (-.55,1.5); 

\draw[thick] (-.35,1) -- (-.35,1.65); 

\draw[thick] (1.5,1) to[out=90,in=0] (-.15,1.5);

\node at (-2,0) {$\scalebox{1}{$=$}$};
\node at (-2,-5.5) {$\scalebox{1}{$\;$}$};
\node at (-2,.5) {$\scalebox{1}{$\;$}$};

\node at (0,-1.2) {$\;$};

\end{tikzpicture}
\]
\end{example} 
The Thompson group is just as good as the braid groups at producing links.
More precisely, Jones proved a result analogous to that of the Alexander theorem.
\begin{theorem}\cite{Jo14}
Given an unoriented link $L$, there exists an element $g$ in $F$ such that $\CL(g)$ is $L$.
\end{theorem}

In analogy to the braid index, it is possible to define a Thompson index.
The $F$-index of a link $L$ is the smallest number of leaves 
 required by each binary tree in a binary tree diagram
 such that $L$ is realised as $\CL(T_+,T_-)$. 
 The $F_3$-index is defined as the smallest number of trivalent vertices plus one
 required by each ternary tree in a binary tree diagram
 such that $L$ is realised as $\CL(T_+,T_-)$. 
Note that the $F_3$-index is defined in terms of trivalent instead of leaves to make it compatible with the $F$-index. In fact, for every binary tree diagram $(T_+,T_-)$ 
the number of trivalent vertices plus $1$ in each ternary tree of $\iota(T_+,T_-)$ is equal to the number of leaves of $T_+$.
The following interesting result was discovered by Golan and Sapir.
\begin{theorem}\label{theoGSindex}\cite{GS}
The $F$-index of a link containing $u$ unlinked unknots and represented by a link diagram with $n$ crossings does not exceed $12n + u + 3$.
\end{theorem}

\section{Positive Thompson knots}\label{sec3}
The positive Thompson knots are those produced by the elements of the monoid of positive words $F_+$.
As explained in Section \ref{sec1}, each of these elements admits a representative whose bottom tree  and  the corresponding graph $\Gamma_-(T_-)$  have the following form 
\begin{eqnarray}\label{Tmeno}
&
\begin{tikzpicture}[x=.55cm, y=.55cm,
    every edge/.style={
        draw,
      postaction={decorate,
                    decoration={markings}
                   }
        }
]
\node (bbb) at (-2,-2) {$\scalebox{1}{$T_-=$}$}; 

\draw[thick] (0,0)--(4,-4)--(8,0);
\draw[thick] (4.5,-3.5)--(1,0);
\draw[thick] (5,-3)--(2,0);
\draw[thick] (7,-1)--(6,0);
\draw[thick] (7.5,-.5)--(7,0);
\draw[thick] (4,-4)--(4,-4.5);

\node (aaaa) at (5,-1) {$\scalebox{1}{$\ldots$}$}; 

\end{tikzpicture} 
\qquad
\begin{tikzpicture}[x=.75cm, y=.75cm,
    every edge/.style={
        draw,
      postaction={decorate,
                    decoration={markings}
                   }
        }
]
\node (aaaa) at (-1.5,0) {$\scalebox{1}{$\Gamma_-(T_-)=$}$}; 
\node (aaaa) at (3,0) {$\scalebox{1}{$\ldots$}$}; 
\node (aaaa) at (3,-1.75) {$\scalebox{1}{$\;$}$}; 

\fill (0,0)  circle[radius=1.5pt];
\fill (1,0)  circle[radius=1.5pt];
\fill (2,0)  circle[radius=1.5pt];
\fill (4,0)  circle[radius=1.5pt];
\fill (5,0)  circle[radius=1.5pt];

\draw[thick] (0,0) to[out=-90,in=-90] (1,0);
\draw[thick] (1,0) to[out=-90,in=-90] (2,0);
\draw[thick] (4,0) to[out=-90,in=-90] (5,0);
\end{tikzpicture}
\end{eqnarray} 
As the form  of the bottom tree is essentially always of the same form
(it depends only on the number of leaves in the upper tree), sometimes we will use the notation $\CL(T_+)$, instead of $\CL(T_+,T_-)$.
Positive Thompson links were the object of study of \cite{AB2}, but before stating the main result of this investigation we recall some preliminary definitions.

\begin{figure}
\phantom{This text will be invisible} 
\[
 \begin{tikzpicture}[x=.35cm, y=.35cm,
    every edge/.style={
        draw,
      postaction={decorate,
                    decoration={markings}
                   }
        }
]

\node at (-4.5,.5) {$\scalebox{1}{$\star 2$}$};


\draw[thick] (0,0)--(1,1);
\draw[thick] (1,0)--(.6,.4);
\draw[thick] (.4,.6)--(0,1);

\draw[thick] (1,0)--(2,1);
\draw[thick] (2,0)--(1.6,.4);
\draw[thick] (1.4,.6)--(1,1);

\draw[thick] (0,1) to[out=90,in=90] (2,1);
\draw[thick] (-1,1) to[out=90,in=90] (3,1);

 \draw[thick] (-1,1) to[out=-90,in=180] (0,0);
 \draw[thick] (3,1) to[out=-90,in=0] (2,0);

\end{tikzpicture} \qquad
\qquad
 \begin{tikzpicture}[x=.35cm, y=.35cm,
    every edge/.style={
        draw,
      postaction={decorate,
                    decoration={markings}
                   }
        }
]

\node at (-6.5,.5) {$\scalebox{1}{$\star$}$};
\node at (-6.5,1.25) {$\scalebox{1}{$1$}$};
\draw[thick] (-6.25,0.5)--(-4,0.5);
\fill (-4,.5)  circle[radius=1.5pt];
\node at (-4,1.25) {$\scalebox{1}{$-2$}$};


\draw[thick] (0,0)--(1,1);
\draw[thick] (1,0)--(.6,.4);
\draw[thick] (.4,.6)--(0,1);

\draw[thick] (2,0)--(2.4,.4);
\draw[thick] (2.6,.6)--(3,1);
\draw[thick] (2,1)--(3,0);
\draw[thick] (2,-1)--(2.4,-.6);
\draw[thick] (2.6,-.4)--(3,0);
\draw[thick] (2,0)--(3,-1);


\draw[thick] (0,1) to[out=90,in=90] (3,1);
\draw[thick] (-1,1) to[out=90,in=90] (4,1);
\draw[thick] (1,1)--(2,1);
\draw[thick] (1,0) to[out=-90,in=180] (2,-1);
\draw[thick] (0,0) to[out=180,in=-90] (-1,1);
\draw[thick] (3,-1) to[out=0,in=-90] (4,1);

\end{tikzpicture} 
\]
\caption{Two arborescent links: one associated with a tree consisting just of the root, another associated with a tree with a root and one leaf.}\label{arbor}
\end{figure}
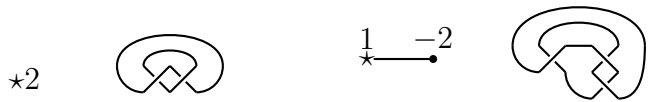

Arborescent tangles are the minimal class of tangles closed under tangle composition, and containing all rational tangles, \cite{Conway}. 
  The closure of an arborescent tangle is described by a finite rooted plane tree with integer vertex weights. Each weight $w$ gives rise to a twist region with $|w|$ crossings. The orientation of these crossings, as well as the interconnections between these twist regions, are determined by the plane tree in the following way. The root vertex corresponds to a horizontal twist region, in which crossings are called positive if their strand going from the bottom left to the top right is above the other strand. 
If the weight is zero, then we have just two horizontal lines.
The vertices adjacent to the root vertex correspond to vertical twist regions attached to this horizontal twist region. The order in which they are attached is determined by the plane cyclic arrangement of the branches around the root vertex. We keep the convention that the overcrossing strand of a positive crossing is going from the bottom left to the top right. In the end, this means that arborescent tangles whose weights carry the same sign give rise to alternating links. The vertices at distance two from the root give again rise to horizontal twist regions, and so on. Two examples are provided in Figure \ref{arbor}. 
For more details and resultst we refer to \cite{BS,Gabai}.
Finally we call a finite rooted plane tree bipartite
 if its vertices have weights $\pm 1$, with the root and all the leaves carrying weight $-1$; all the vertices of weight 1 have degree $2$; there are no edges between vertices with the same weight.

We are now in a position to state the main result of \cite{AB2}, which is kind of an Alexander theorem for $F_+$.
\begin{theorem}\cite{AB2}
The set of positive Thompson links coincides with the set of closures of bipartite arborescent tangles.
\end{theorem}

In general links produced from $F$ and $F_3$ do not posses a natural orientation, for this reason Jones introduced the so-called oriented subgroups.
We will define them in the next section.

 \section{The oriented subgroups $\vec{F}$ and $\vec{F}_3$}\label{sec4}
In this section we introduce two new subgroups, one of $F$ and one of $F_3$.
They are interesting on their own, but we will need them to produce oriented links.
 The (binary) oriented subgroup $\vec{F}=\vec{F}_2$ is a subgroup of $F$, while $\vec{F}_3$ is a subgroup of $F_3$.

 First we describe an alternative method to obtain the Tait graph associated with each ternary tree diagram.
We will ignore the sign of the edges as it is not relevant for defining the oriented subgroups.
Given a ternary tree diagram $(T_+,T_-)$, we call this graph $\Gamma(T_+,T_-)$ the planar graph of $(T_+,T_-)$.
 We imagine $(T_+,T_-)$ sitting in the strip bounded by the lines $y=1$ and $y=-1$.
This strip is $2$-colourable. We use two colors: black and white, the left-most region is black.
The vertices of  $\Gamma(T_+,T_-)$  sit on the $x$-axis, precisely 
on $-1/2+2\IN_0 :=\{-1/2,1+1/2, 3+1/2, \ldots \}$ 
and there is precisely one vertex for  every black region. 
We draw an edge between two black regions whenever they meet at a $4$-valent vertex. 
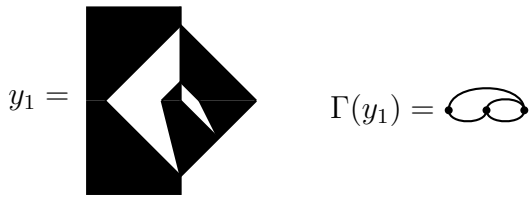
\[
\begin{tikzpicture}[x=.5cm, y=.5cm,
    every edge/.style={
        draw,
      postaction={decorate,
                    decoration={markings}
                   }
        }
]

\node at (-.75,0) {$\scalebox{1}{$y_1=$}$};

\draw[thick] (1,0) -- (3,2)--(5,0);
 \draw[thick] (3,0) -- (3,2);
\draw[thick] (2.5,0) -- (3,.5)--(3.5,0);

\fill [color=black,opacity=0.2]   (.5,0)--(.5,2.5)--(3,2.5)--(3,2)--(1,0); 
\fill [color=black,opacity=0.2]   (2.5,0)--(3,.5)--(3,0); 
\fill [color=black,opacity=0.2]   (3,2)--(5,0)--(3.5,0)--(3,0.5)--(3,2); 

\draw[thick] (3,-2)--(3,-2.5);
\draw[thick] (3,2)--(3,2.5);

\draw[thick] (1,0) -- (3,-2)--(5,0);
\draw[thick] (3,0) -- (4,-1);
\draw[thick] (3,-2) -- (2.5,0);
\draw[thick] (4,-1) -- (3.5,0);

\fill [color=black,opacity=0.2]   (.5,0)--(.5,-2.5)--(3,-2.5)--(3,-2)--(1,0); 
\fill [color=black,opacity=0.2]   (2.5,0)--(3,-2)--(4,-1)--(3,0); 
\fill [color=black,opacity=0.2]   (5,0)--(4,-1)--(3.5,0)--(5,0);

\end{tikzpicture}
\qquad
\begin{tikzpicture}[x=.5cm, y=.5cm,
    every edge/.style={
        draw,
      postaction={decorate,
                    decoration={markings}
                   }
        }
]

\draw[thick] (1,0) to[out=90,in=90] (3,0);
\draw[thick] (3,0) to[out=90,in=90] (2,0);

\draw[thick] (2,0) to[out=-90,in=-90] (3,0);
\draw[thick] (1,0) to[out=-90,in=-90] (2,0);

\fill (1,0)  circle[radius=1.5pt];
\fill (2,0)  circle[radius=1.5pt];
\fill (3,0)  circle[radius=1.5pt];
  
\node at (-.75,0) {$\scalebox{1}{$\Gamma(y_1)=$}$};
\node at (-1.25,-2) {$\scalebox{1}{$\;$}$};

\end{tikzpicture}
\]
The binary and the ternary oriented subgroup $\vec{F}_3$ can be defined as
\begin{align}
\vec{F}&=\vec{F}_2:=\{(T_+,T_-)\in F\; | \; \Gamma(\iota(T_+,T_-)) \text{ is $2$-colorable}\}\\
\vec{F}_3&:=\{(T_+,T_-)\in F_3\; | \; \Gamma(T_+,T_-) \text{ is $2$-colorable}\}
\end{align}
\begin{convention}
We denote the colors used for the vertices of $\Gamma(T_+,T_-)$ by
 $+$ and $-$.
The graph $\Gamma(T_+,T_-)$ is always connected and, therefore, if it is $2$-colorable there are only $2$ possible colorings: one where the leftmost vertex has color $+$, one with color $-$. By convention we choose always choose the first of these colorings.
\end{convention}
These groups were introduced by Jones in 2014 \cite{Jo14} and 2018 \cite{Jo18}, respectively.
The (binary) oriented subgroup was first studied by Golan and Sapir in \cite{GS, GS2}, who determined its generators
  $x_0x_1$, $x_1x_2$, $x_2x_3$.
Moreover, they also discovered that the map induced by 
\begin{align*}
\alpha: & \; \vec{F}\to F_3\\
&x_0x_1\mapsto y_0\\
&x_1x_2\mapsto y_1\\
&x_2x_3\mapsto y_2
\end{align*}
is an isomorphism. A pictorial interpretation of this isomorphism was later found by Ren in \cite{Ren}, where he realised that this map can be obtained by taking a ternary tree diagram and replacing each trivalent vertex with a suitable tree with $3$ leaves, see Figure \ref{fig-ren-map}.

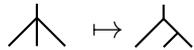
\begin{figure}
\phantom{This text will be invisible} 
\[\begin{tikzpicture}[x=.75cm, y=.75cm,
    every edge/.style={
        draw,
      postaction={decorate,
                    decoration={markings}
                   }
        }
]

\draw[thick] (0,0)--(.5,.5)--(1,0);
\draw[thick] (0.5,0.75)--(.5,0);
\node at (1.75,0.25) {$\scalebox{1}{$\mapsto$}$};

\end{tikzpicture}
\begin{tikzpicture}[x=.75cm, y=.75cm,
    every edge/.style={
        draw,
      postaction={decorate,
                    decoration={markings}
                   }
        }
]

\draw[thick] (0.5,0.75)--(.5,.5);
\draw[thick] (0.5,0)--(.75,.25);
\draw[thick] (0,0)--(.5,.5)--(1,0);
\end{tikzpicture}
\]
 \caption{Ren's map 
 for $\vec{F}$.} 
  \label{fig-ren-map}
\end{figure} 
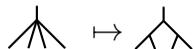
\begin{figure}
\[\begin{tikzpicture}[x=.75cm, y=.75cm,
    every edge/.style={
        draw,
      postaction={decorate,
                    decoration={markings}
                   }
        }
]

\draw[thick] (0,0)--(.5,.5)--(1,0);
\draw[thick] (0.5,0.5)--(.5,.75);
\draw[thick] (0.5,0.5)--(.35,0);
\draw[thick] (0.5,0.5)--(.65,0);
 \node at (1.75,0.25) {$\scalebox{1}{$\mapsto$}$};

\end{tikzpicture}
\begin{tikzpicture}[x=.75cm, y=.75cm,
    every edge/.style={
        draw,
      postaction={decorate,
                    decoration={markings}
                   }
        }
]

\draw[thick] (0.5,0.75)--(.5,.5);
\draw[thick] (0.65,0)--(.75,.25);
\draw[thick] (0.35,0)--(.25,.25);
\draw[thick] (0,0)--(.5,.5)--(1,0);
 \end{tikzpicture}
\]
\caption{The map 
 for  $\CF$.}
  \label{fig-ren-map-2}
\end{figure} 

One may define a \emph{weight} $\omega$ set of finite binary words $W_2$ with values in $\IZ_2$ by the formula $\omega(a_1\ldots a_n):=\sum_{i=0}^n a_i$
and a subset of dyadic rationals
$$
S:=\{t\in W_2 \; | \;  \omega(t)=0\} .
$$
\begin{theorem}\cite{GS}
The oriented subgroup $\vec{F}$ is the stabiliser subgroup $\stab(S)$.
\end{theorem}
Thanks to this characterisation Golan and Sapir  were able to prove that $\vec{F}$ coincided with its commensurator
 and thus the corresponding quasi-regular representation is 
irreducible.

The oriented subgroup $F$ is also interesting because it gave rise to a novel example of maximal subgroup of infinite index in $F$.
Before this group was defined, the only known subgroups of this type were the so-called parabolic subgroups, that is the stabilisers of points $\stab(t)$, $t\in (0,1)$,
under the natural action of $F$ on $(0,1)$ that were studied in \cite{Sav, Sav2}.  
Golan and Sapir in \cite{GS2} proved that $\vec{F}$ sits inside the rectangular subgroup 
$K_{(1,2)}$ and it is maximal and of infinite index in it. 
By exhibiting the explicit isomorphism $\beta: K_{(1,2)}\to F$ induced by $\beta(x_0x_2)=x_0$, $\beta(x_1x_2)=x_1$, they were able to show that 
$\beta(\vec{F})$ was a maximal subgroup of infinite index in $F$ distinct from the parabolic subgroups.
%

Once the ternary oriented subgroup was introduced, we extended some of  these results to it.
More precisely, we found a set of  generators for $\vec{F}_3$.
\begin{theorem}\cite{TV}
The ternary oriented Thompson group $\vec{F}_3$ is generated by the following elements
\begin{align*}
& y_{2i+1}^2, y_{2i}y_{2i+2}, y_{2i}y_{2i+3}
\qquad i = 0, 1, 2.
\end{align*}
 \end{theorem}
\begin{question}
Is $\vec{F}_3$ isomorphic to a Brown-Thompson group or other known groups?
\end{question}
One might also ask the following question (which was originally asked in \cite[Problem 1]{TV}).
\begin{question}
Is $\vec{F}_3$ finitely presented?
\end{question}
Then we also realised that $\vec{F}_3$ is the stabiliser of a suitable subset. 
For this we needed a new weight on the set of ternary words.

Consider a tree in the upper half-plane and its leaves on the $x$-axis as usual.
To each vertex 
 $v$ of a tree we associate a natural number $c(v)$ which we call its {\bf weight}, as follows. 
Given a vertex, there exists a unique minimal path from the root of the tree to the vertex. This path is made by a collection of left, middle, right edges, and may be represented by a word $w_1 1w_2 1\cdots 1 w_n$ in the letters
$\{0, 1, 2\}$ ($0$ stands for a left edge, $1$ for a middle edge, $2$ for a right edge), where $w_1, \ldots , w_{n-1}$ are words that do not contain the letter $1$, $w_n$ can have $1$ only as its last letter.
We call $\{w_{2k+1}\}_{k\geq 0}$ the odd words and $\{w_{2k}\}_{k\geq 0}$ the even words.
The weight of $v$ 
is the sum of the number of digits equal to $1$, plus the number of digits equal to $2$ in the odd words, plus the number of digits equal to $0$ in the even words.
When we compute the weight of a leaf in a tree diagram, 
sometimes we use the symbol $c_+$ or $c_-$ to distinguish which tree we are considering ($c_+$ for the top tree, $c_-$ for the reflected bottom tree).
 
\begin{theorem}\cite{TV}
The  ternary oriented Thompson group $\vec{F}_3$ is the stabilizer of the following subset of the triadic fractions
$$
Z:=\{ .a_1a_2\cdots a_n \; | \; \# \textrm{ of $1$'s is even}, c(.a_1a_2\cdots a_n) \textrm{ is even}\} \; .
$$
\end{theorem}
 
\begin{figure}
\[
\begin{tikzpicture}[x=.6cm, y=.6cm, every path/.style={
 thick}, 
every node/.style={transform shape, knot crossing, inner sep=1.5pt}]

\node (aaaa) at (-1.5,1) {$\scalebox{1}{$+=$}$}; 

\node (a1) at (0,0) {};
\node (a2) at (2,2) {};

\node (b3) at (0,2) {};
\node (b2) at (1,1) {};
\node (b1) at (2,0) {};

\draw[->] (a1.center) .. controls (a1.4 north east) and (a2.4 south west) .. (a2.center);
\draw (b1.center) .. controls (b1.4 north west) and (b2.4 south east) .. (b2);
\draw[->] (b2) .. controls (b2.4 north west) and (b3.4 south east) .. (b3.center);

\end{tikzpicture}
\qquad\qquad\qquad
\begin{tikzpicture}[x=.6cm, y=.6cm, every path/.style={
 thick}, 
every node/.style={transform shape, knot crossing, inner sep=1.5pt}]

\node (aaaa) at (-1.5,1) {$\scalebox{1}{$-=$}$}; 

\node (a1) at (0,0) {};
\node (a2) at (2,2) {};

\node (b3) at (0,2) {};
\node (b2) at (1,1) {};
\node (b1) at (2,0) {};

\draw[->] (b1.center) .. controls (b1.4 north west) and (b3.4 south east) .. (b3.center);
\draw (a1.center) .. controls (a1.4 north east) and (b2.4 south west) .. (b2);
\draw[->] (b2) .. controls (b2.4 north east) and (a2.4 south west) .. (a2.center);

\end{tikzpicture}
\]
\caption{A positive and a negative crossing in an oriented link.} \label{figpos}
\end{figure}
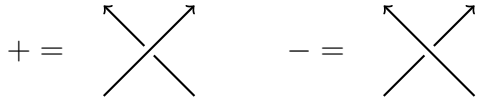

\section{The construction of knots: the oriented case}\label{sec5}
Links produced from the oriented subgroups admit a natural orientation. Recall that the Tait diagrams $\Gamma(T_+,T_-)$ of elements in $\vec{F}$ are  $2$-colorable
(the colors being $\{+, -\}$ and the left-most vertex having colour $+$).

Given $(T_+,T_-)$, if we shade  the   link diagram $\CL(T_+,T_-)$ in black and white (we adopt the convention that the colour of the unbounded region is white),
by construction the vertices of the graph $\Gamma(T_+, T_-)$ sit in the black regions and each one has been assigned with a colour $+$ or $-$. 
These colours determine an orientation of the surface and of the boundary ($+$ means that the region is positively oriented). 
It can be easily seen that the graph $\Gamma(\iota(x_0x_1))$ is $2$-colorable and thus $x_0x_1$ is in $\vec{F}$. 
Here is the oriented link associated with $x_0x_1$.
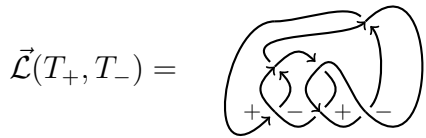
\[
\begin{tikzpicture}[x=.6cm, y=.6cm, every path/.style={
 thick}, 
every node/.style={transform shape, knot crossing, inner sep=1.5pt}]

\node (aaaa) at (-3,1) {$\scalebox{1}{$\vec\CL(T_+,T_-)=$}$}; 

\node (a2) at (1,0) {};
\node (a3) at (2,0) {};
\node (a4) at (3,0) {};
\node (b2) at (1,1) {};
\node (b3) at (2,1) {};
\node (b4) at (3,2) {};

\node (appo) at (1,2) {};

\draw[->] (a2.center) .. controls (a2.4 north west) and (b2.4 south west) .. (b2.center);
\draw[->] (a2) .. controls (a2.4 north east) and (b2.4 south east) .. (b2);

\draw[->] (b2) .. controls (b2.8 north west) and (b4.4 south west) .. (b4.center);

\draw[->] (b2.center) .. controls (b2.4 north east) and (b3.4 north west) .. (b3);
\draw (a2.center) .. controls (a2.4 south east) and (a3.4 south west) .. (a3);
\draw (a3) .. controls (a3.4 north east) and (b3.4 south east) .. (b3);
\draw[<-] (a3.center) .. controls (a3.4 north west) and (b3.4 south west) .. (b3.center);
\draw[<-] (a2) .. controls (a2.16 south west) and (appo.8 west) .. (appo.center);
\draw (appo.center) .. controls (appo.4 east) and (b4.4 north west) .. (b4);

\draw (a3.center) .. controls (a3.4 south east) and (a4.4 south west) .. (a4);
\draw (a4.center) .. controls (a4.16 south east) and (b4.16 north east) .. (b4.center);
\draw[->] (a4) .. controls (a4.4 north east) and (b4.4 south east) .. (b4);

\draw (b3.center) .. controls (b3.4 north east) and (a4.4 north west) .. (a4.center);

 \node (x2) at (.5,0) {$\scalebox{0.75}{$+$}$};
\node (x3) at (1.45,0) {$\scalebox{0.75}{$-$}$};
\node (x4) at (2.5,0) {$\scalebox{0.75}{$+$}$};
\node (x5) at (3.4,0) {$\scalebox{0.75}{$-$}$};
 
\end{tikzpicture} 
\]
Also in the setting of the oriented subgroups, Jones proved a result analogous to Alexander's theorem.
\begin{theorem}\cite{Jo14}
Given an oriented link $\vec{L}$, there exists an element $g$ in $\vec{F}$ such that $\vec{\CL}(g)$ is $\vec{L}$ up to disjoint union with unknots.
\end{theorem}
In \cite{Jo18} Jones asked whether the previous theorem could be improved and each oriented link could be exactly reproduced. An answer was provided in \cite{A}.
\begin{theorem}\cite{A}
Given an oriented link $\vec{L}$, there exists an element $g$ in $\vec{F}$ such that $\vec{\CL}(g)$ is $\vec{L}$.
\end{theorem}

There is an interesting monoid in $\vec{F}$, namely the monoid $\vec{F}_+:=\vec{F}\cap F_+$.
It is well known that every element of the braid group may be expressed as the product of a positive braid and the inverse of a positive braid.
 Similarly, for the oriented Thompson group we have the following result which is due to Ren, \cite{Ren}.
\begin{proposition}\cite{Ren}
For every $g\in\vec{F}$, there exist $g_+$, $g_-\in\vec{F}_+$ such that 
$g=g_+(g_-)^{-1}$.
\end{proposition}

We call the links produced by the monoid $\vec{F}_+$ the positive oriented Thompson links. 
Recall that an oriented link is called positive if it admits a link diagram where all it crossings are positive in the sense of Figure \ref{figpos}.
\begin{theorem}\cite{AB1}
The positive oriented Thompson links are positive.
\end{theorem}
As in for $F$ and $F_3$, it is possible to define the $\vec{F}$ and $\vec{F}_3$ indices.
\begin{question}
Does Theorem \ref{theoGSindex} extend to the $\vec{F}$ and $\vec{F}_3$ indices?
\end{question} 

\section{The $3$-colorable subgroup $\CF$} \label{sec6}
 Another subgroup introduced by Jones is the so-called  $3$-colorable subgroup $\CF$.
As before, any binary tree diagram partitions the strip bounded by the lines $y=1$ and $y=-1$ in regions.
This strip may or may not be $3$-colorable, i.e., it may or may not be 
 possible to assign the colors $\IZ_3=\{0,1, 2\}$ to the regions of the strip in such a way that if two regions share an edge, they   have different colors.  \begin{convention}
 If the strip is $3$-colorable, we adopt this convention: we assign the following colors to the regions near the roots 
  \begin{eqnarray}\label{convention-colors}
\begin{tikzpicture}[x=.5cm, y=.5cm,
    every edge/.style={
        draw,
      postaction={decorate,
                    decoration={markings}
                   }
        }
] 

 \draw[thick] (0,0) -- (1,1)--(2,0);
 \draw[thick] (1,1) -- (1,2);



\node at (0,1) {$0$};
\node at (2,1) {$1$};
\node at (1,0.2) {$2$};

 \draw[thick] (-1,2) -- (3,2);


\end{tikzpicture}
\qquad\qquad 
\begin{tikzpicture}[x=.5cm, y=.5cm,
    every edge/.style={
        draw,
      postaction={decorate,
                    decoration={markings}
                   }
        }
] 

 \draw[thick] (-1,0) -- (3,0);

 \draw[thick] (0,2) -- (1,1)--(2,2);
 \draw[thick] (1,1) -- (1,0);



\node at (0,1) {$0$};
\node at (2,1) {$1$};
\node at (1,1.8) {$2$};

\end{tikzpicture}
\end{eqnarray}
Once we make this convention, if the strip is $3$-colourable, 
there exists a unique colouring. 
\end{convention}
The \textbf{$3$-colorable subgroup} $\CF$  consists of the elements of $F$ whose corresponding  strip is $3$-colorable. 
For example, this is the strip corresponding to $w_0:=x_0^2x_1x_2^{-1}$ (which is $3$-colorable)
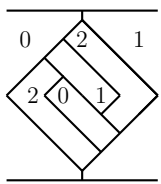
\[
\begin{tikzpicture}[x=.5cm, y=.5cm,
    every edge/.style={
        draw,
      postaction={decorate,
                    decoration={markings}
                   }
        }
]

\draw[thick] (0,0)--(2,2)--(4,0);
\draw[thick] (0,0)--(2,-2)--(4,0);
\draw[thick] (1,0)--(1.5,0.5)--(2,0);
\draw[thick] (1.5,0.5)--(1,1);
\draw[thick] (3,0)--(1.5,1.5);
\draw[thick] (3,0)--(2.5,-.5); 
\draw[thick] (4,0)--(2,2);
\draw[thick] (4,0)--(3.5,-.5);
\draw[thick] (1,0)--(2.5,-1.5);
\draw[thick] (2,0)--(3,-1);

\draw[thick] (2,2)--(2,2.25);
\draw[thick] (2,-2)--(2,-2.25);

\node at (0,-1.2) {$\;$};
 \draw[thick] (0,2.25)--(4,2.25);

 \draw[thick] (0,-2.25)--(4,-2.25);


\node at (0.5,1.5) {$\scalebox{.75}{$0$}$};
\node at (2,1.5) {$\scalebox{.75}{$2$}$};
\node at (3.5,1.5) {$\scalebox{.75}{$1$}$};
\node at (0.7,0) {$\scalebox{.75}{$2$}$};
\node at (1.5,0) {$\scalebox{.75}{$0$}$};
\node at (2.5,0) {$\scalebox{.75}{$1$}$};

\end{tikzpicture}
\]
The $3$-colorable subgroup actually provides a copy of a Borwn-Thompson group inside $F$.
\begin{theorem} \cite{Ren}
The subgroup  $\CF$
is isomorphic to the Brown-Thompson group $F_4$   by means of the isomorphism (with domain $F_4$ and range $\CF$) obtained
 by  replacing every $5$-valent vertex of $4$-ary trees by the complete binary tree with $4$ leaves 
(see Figure \ref{fig-ren-map-2}).
The images of   $y_0$, $y_1$, $y_2$, $y_3$ yield the following elements
$w_0:=x_0^2x_1x_2^{-1}$,
$w_1:=x_0x_1^2x_0^{-1}$,
$w_2:=x_1^2x_3x_2^{-1}$,
$w_3:=x_2^2x_3x_4^{-1}$ (see Figure \ref{genCF}). 
\end{theorem} 

\begin{figure}
\[\begin{tikzpicture}[x=.75cm, y=.75cm,
    every edge/.style={
        draw,
      postaction={decorate,
                    decoration={markings}
                   }
        }
]

\draw[thick] (0,0)--(.5,.5)--(1,0);
\draw[thick] (0.5,0.5)--(.5,.75);
\draw[thick] (0.5,0.5)--(.35,0);
\draw[thick] (0.5,0.5)--(.65,0);
\node at (0,-1.2) {$\;$};
\node at (1.75,0.25) {$\scalebox{1}{$\mapsto$}$};

\end{tikzpicture}
\begin{tikzpicture}[x=.75cm, y=.75cm,
    every edge/.style={
        draw,
      postaction={decorate,
                    decoration={markings}
                   }
        }
]

\draw[thick] (0.5,0.75)--(.5,.5);
\draw[thick] (0.65,0)--(.75,.25);
\draw[thick] (0.35,0)--(.25,.25);
\draw[thick] (0,0)--(.5,.5)--(1,0);
\node at (0,-1.2) {$\;$};
\end{tikzpicture}
\]
\caption{Ren's map  $\Phi$ 
 from the set of $4$-ary trees to binary trees.  }
  \label{fig-ren-map-2}
\end{figure}
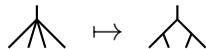  
  \begin{figure}
\phantom{This text will be invisible} 
\[
\begin{tikzpicture}[x=.75cm, y=.75cm,
    every edge/.style={
        draw,
      postaction={decorate,
                    decoration={markings}
                   }
        }
]

\node at (-1.5,-.25) {$\scalebox{1}{$w_0=$}$};

\draw[thick] (0,0)--(2,2)--(4,0);
\draw[thick] (0,0)--(2,-2)--(4,0);
\draw[thick] (1,0)--(1.5,0.5)--(2,0);
\draw[thick] (1.5,0.5)--(1,1);
\draw[thick] (3,0)--(1.5,1.5);
\draw[thick] (3,0)--(2.5,-.5); 
\draw[thick] (4,0)--(2,2);
\draw[thick] (4,0)--(3.5,-.5);
\draw[thick] (1,0)--(2.5,-1.5);
\draw[thick] (2,0)--(3,-1);

\draw[thick] (2,2)--(2,2.25);
\draw[thick] (2,-2)--(2,-2.25);

\node at (0,-1.2) {$\;$};
\end{tikzpicture}
\qquad 
\begin{tikzpicture}[x=.75cm, y=.75cm,
    every edge/.style={
        draw,
      postaction={decorate,
                    decoration={markings}
                   }
        }
]

\node at (-1.5,-.25) {$\scalebox{1}{$w_1=$}$};

\draw[thick] (0,0)--(2,2)--(4,0);
\draw[thick] (0,0)--(2,-2)--(4,0);

\draw[thick] (1,0)--(1.5,0.5)--(2,0);
\draw[thick] (2,0)--(3,-1);
\draw[thick] (1,0)--(.5,-.5);
\draw[thick] (1.5,.5)--(2,1)--(3,0);
\draw[thick] (3,0)--(3.5,-.5);
\draw[thick] (2,1)--(1.5,1.5);

\draw[thick] (2,2)--(2,2.25);
\draw[thick] (2,-2)--(2,-2.25);

\node at (0,-1.2) {$\;$};
\end{tikzpicture}
\]

\[
\begin{tikzpicture}[x=.75cm, y=.75cm,
    every edge/.style={
        draw,
      postaction={decorate,
                    decoration={markings}
                   }
        }
]

\node at (-2.5,-.25) {$\scalebox{1}{$w_2=$}$};
 
\draw[thick] (-1,0)--(1.5,2.5)--(2,2);
\draw[thick] (-1,0)--(1.5,-2.5)--(2,-2);

\draw[thick] (0,0)--(2,2)--(4,0);
\draw[thick] (0,0)--(2,-2)--(4,0);

\draw[thick] (2,0)--(2.5,.5)--(3,0);

\draw[thick] (1,0)--(.5,.5);
\draw[thick] (2.5,.5)--(1.5,1.5);

\draw[thick] (3,0)--(3.5,-.5);
\draw[thick] (2,0)--(1.5,-.5)--(1,0);
\draw[thick] (1.5,-.5)--(2.5,-1.5);

\draw[thick] (1.5,2.5)--(1.5,2.75);
\draw[thick] (1.5,-2.5)--(1.5,-2.75);

\node at (0,-1.2) {$\;$};
\end{tikzpicture}
\qquad 
\begin{tikzpicture}[x=.75cm, y=.75cm,
    every edge/.style={
        draw,
      postaction={decorate,
                    decoration={markings}
                   }
        }
]

\node at (-3.5,-.25) {$\scalebox{1}{$w_3=$}$};
 
\draw[thick] (-1,0)--(1.5,2.5)--(2,2);
\draw[thick] (-1,0)--(1.5,-2.5)--(2,-2);

\draw[thick] (-2,0)--(1,3)--(1.5,2.5);
\draw[thick] (-2,0)--(1,-3)--(1.5,-2.5);

\draw[thick] (0,0)--(2,2)--(4,0);
\draw[thick] (0,0)--(2,-2)--(4,0);
\draw[thick] (1,0)--(1.5,0.5)--(2,0);
\draw[thick] (1.5,0.5)--(1,1);
\draw[thick] (3,0)--(1.5,1.5);
\draw[thick] (3,0)--(2.5,-.5); 
\draw[thick] (4,0)--(2,2);
\draw[thick] (4,0)--(3.5,-.5);
\draw[thick] (1,0)--(2.5,-1.5);
\draw[thick] (2,0)--(3,-1);

\draw[thick] (1,3)--(1,3.25);
\draw[thick] (1,-3)--(1,-3.25);

\node at (0,-1.2) {$\;$};
\end{tikzpicture}
\]
 \caption{The generators of $\CF$. 
 }\label{genCF}
\end{figure}
The 3-colorable subgroup is the intersection of the stabilizers of three subsets of dyadic rationals.
\begin{theorem}\cite{TV2}
For a binary word $a_1a_2\ldots a_n$  set
$$
\omega(a_1a_2\ldots a_n):= \sum_{i=1}^n (-1)^ia_i \in\IZ_3\, .$$
$$
S_i:=\{t\in (0,1)\cap \IZ[1/2]\; | \; \omega(t)= i\} \qquad i\in \IZ_3\; 
$$
where 
$\equiv_3$ is the equivalence modulo $3$. 
Then it holds 
$$
\mathcal{F}=\cap_{i \in \mathbb{Z}_3} {\rm Stab}(S_i)\, .
$$
\end{theorem}
Simple computations show that the $3$-colorable subgroup $\CF$ is contained in the rectangular subgroup $K_{(2,2)}\cong F$, but unlike $\vec{F}$ it is not maximal in it. However, there are still three maximal subgroups of $K_{(2,2)}$ of infinite index, namely $M_0:=\langle x_0^2, \CF\rangle$,
$M_1:=\langle x_1^2, \CF\rangle$, and $M_2:=\langle \sigma(x_1)^2, \CF\rangle$, \cite{TV2}.

  Connections between  $\CF$ and Jones's construction of knots have not been explored yet, nevertheless it is natural to ask the following question.
 \begin{question}
Do the elements of  $\CF$ produce all unoriented knots and links?
\end{question}

\section*{Acknowledgements}
The author is grateful to Stefano Rossi and the anonymous referee for their valuable comments on the first version of this paper.

\end{document}